\documentclass[a4paper,12pt]{article} 

\usepackage{color}

\usepackage{amsfonts, amsmath, amsthm, amssymb}
\usepackage[T1]{fontenc}
\usepackage[cp1250]{inputenc}
\usepackage{xcolor}
\usepackage{graphicx}
\usepackage{amssymb}
\usepackage{amsmath}
\usepackage{mathptmx}
\usepackage{helvet}
\usepackage{courier}
\usepackage{txfonts}
\usepackage{tikz} 
\usetikzlibrary{arrows}
\usepackage{type1cm}

\usepackage{verbatim}

\usepackage{graphicx}
\usepackage{epsfig,amscd,amssymb,amsxtra,amsmath,amsthm}
\usepackage{type1cm}
\usepackage[T1]{fontenc}
\usepackage{graphics}
\usepackage[mathscr]{eucal}
\usepackage[all]{xy}
\usepackage{amsmath,amscd}

\newcommand{\orbit}{\mathcal{O}^{\oplus}}

\newtheorem{theorem}{Theorem}[section]
\newtheorem{proposition}[theorem]{Proposition}
\newtheorem{definition}[theorem]{Definition}
\newtheorem{lemma}[theorem]{Lemma}

\newtheorem{example}[theorem]{Example}

\newtheorem{problem}[theorem]{Problem}
\newtheorem{observation}[theorem]{Observation}

\newcommand{\diam} {\mathop{\rm diam}\nolimits}

\newcommand{\Cl}  {\mathop{\rm Cl}\nolimits}
\newcommand{\Int}  {\mathop{\rm Int}\nolimits}

\begin{document}

\def\joinrel{\mkern-3mu}
\newcommand{\varproj}{\displaystyle \lim_{\multimapinv\joinrel-\joinrel-}}

\title{Transitive mappings on the Cantor fan}
\author{Iztok Bani\v c, Goran Erceg, Judy Kennedy,  Chris Mouron and Van Nall}
\date{}

\maketitle

\begin{abstract}

{Many  continua that admit a transitive homeomorphism may be found in the literature.  The circle  is probably  the  simplest non-degenerate continuum that admits such a homeomorphism. On the other hand, most of the known examples of such continua have a  complicated topological structure. For example, they are {indecomposable} (such as the pseudo-arc or the Knaster bucket-handle continuum), or they are {not indecomposable} but have some other complicated topological structure, such as  a dense set of ramification points (such as the Sierpi\' nski carpet)  or a dense set of end-points (such as the Lelek fan).  In this paper, we continue our mission of finding continua with simpler topological structures that admit a transitive homeomorphism.} We construct a transitive homeomorphism on the Cantor fan.  

{In our approach, we use four different techniques, each of them giving a unique construction of a transitive homeomorphism on the Cantor fan:} two techniques using quotient spaces of products of compact metric spaces and Cantor sets, and two using Mahavier products of closed relations on compact metric spaces.  {We also demonstrate how our technique using Mahavier products of closed relations may be used to } construct a transitive function $f$ on a Cantor fan $X$ such that $\varprojlim(X,f)$ is a Lelek fan. 
\end{abstract}
\-

\section{Introduction}
In topological dynamical systems theory,  an interesting behavior of a dynamical system $(X,f)$ is often inherited from topological properties of  $X$ and from the  properties of the mapping $f$.  One of the commonly studied properties in the theory of topological dynamical systems is transitivity of dynamical systems $(X,f)$.  {In this paper, we are interested in transitive  topological dynamical systems $(X,f)$, where $X$ is a continuum and $f$ is a homeomorphism. We say that a continuum $X$ admits a transitive homeomorphism, if there is a homeomorphism $f$ on $X$ such that $(X,f)$ is transitive. In the literature, there are many examples of continua that admit a transitive homeomorphism (such as the circle, the pseudo-arc, the Sierpi\' nski carpet, or the Lelek fan). It is often the case that when constructing a non-trivial continuum $X$ that admits a transitive homeomorphism, the continuum $X$ should have a complicated topological structure.  Our main focus in this paper  is to present various techniques that allow us to construct examples of simple continua that admit transitive homeomorphisms.  Then we demonstrate how they can be applied to show that the Cantor fan admits a transitive homeomorphism.  }

We present two different techniques for constructing a transitive dynamical system $(X,f)$ on a compact metric space (or a continuum) $X$:
\begin{itemize}
\item {\bf Our technique using quotient spaces. } For a given compact metric space $X$ and for a given family $\mathcal F$ of continuous functions on $X$,  we construct a transitive dynamical system $(Y,f)$ (depending on $X$ and $\mathcal F$) by defining $Y$ as a quotient space of the topological product $X\times C$, where $C$ is a Cantor set obtained as an infinite product of a finite set. The  transitive mapping $f$ on the space $Y$ is then defined using the family $\mathcal F$ as well as the shift map on $C$.  We use this technique to present two different proofs that there is a transitive homeomorphism on the Cantor fan.  
\item {\bf Our technique using closed relations on compact metric spaces.} When constructing a dynamical system  with closed relations $(X,F)$,  two  standard topological  dynamical systems $(X_F^+,\sigma_F^+)$ and $(X_F,\sigma_F)$ are constructed  at the same time.  This is a new way of constructing topological dynamical systems with interesting properties. Therefore, it is only natural to study the properties of $(X,F)$ that imply interesting topological or dynamical properties of the systems $(X_F^+,\sigma_F^+)$ and $(X_F,\sigma_F)$, where $\sigma_F$ is a homeomorphism while $\sigma_F^+$ is usually not.  Here we show that the shift homeomorphism on $X_F$ is transitive if and only if the shift map on $X_F^+$ is transitive.  We also study closed relations on $X$ that are unions of graphs of continuous functions.  Then, we apply these results to obtain two transitive homeomorphisms on the Cantor fan.

At the end, we also study  $\sigma$-transitivity of dynamical systems  with closed relations.  One of the results  we  prove is that $\sigma$-transitivity of $(X,F)$ is equivalent to the transitivity of the dynamical system $(X_F^+,\sigma_F^+)$.  We conclude the paper by giving an illustrative  example showing that there is a transitive function $f$ on a Cantor fan $X$ such that $\varprojlim(X,f)$ is a Lelek fan. 
\end{itemize} 

As already mentioned,  to construct an interesting  dynamical system $(X,\varphi)$, where $X$ is a continuum and $\varphi$ a homeomorphism,  the continuum $X$ often needs to have some complicated topological structure (see \cite{banic2,barge,jan,jan2,jan3,jan4,HM,chris2,chris3,chris4,chris5,oprocha,seidler}, where { examples of} such continua are obtained).  In this paper, we are interested in the transitivity of such  dynamical systems $(X,\varphi)$.  To date,  many transitive dynamical systems on continua have been constructed.   For example, on a circle, any irrational rotation is a transitive homeomorphism.  
 In \cite{judy}, J. ~Kennedy  and in \cite{minc},  P. ~Minc and W. ~R.~ R. ~Transue constructed independently transitive homeomorphisms on the pseudo-arc.  P. ~Minc and W. ~R.~ R. ~Transue's construction is based on inverse limits of unit intervals.  In \cite{cinc},  J.~\v Cin\v c and P.~Oprocha constructed a {complete space of transitive homeomorphisms with various additional dynamical properties} on the pseudo-arc. This generalizes the  results from \cite{judy} and \cite{minc}. 
In \cite{handel},  M.~Handel constructed transitive homeomorphisms on the pseudo-circle.  In \cite{chris1},  V. ~Mart\' inez-de-la-Vega,  J.~M.~Mart\' inez-Montejano and  C.~Mouron constructed a transitive mapping (actually mixing) on Wa\. zewski's universal dendrite,  which,  after taking the inverse limit produces a mixing homeomorphism on a hereditarily decomposable tree-like continuum by using the shift,    and in \cite{jan4},  J.~ Boro\' nski and P.~Oprocha constructed a transitive homeomorphism on the Sierpi\' nski carpet.  In \cite{banic2}, I.~Bani\v c, G.~Erceg and J.~Kennedy present a transitive homeomorphism on the Lelek fan. Their homeomorphism has non-zero entropy.  In \cite{oprocha},  motivated by \cite{banic2},  P.~Oprocha constructed  a transitive homeomorphism on the Lelek fan; the entropy of his homeomorphism is zero.  

{ In this paper,  we develop various techniques that are applied to study transitive mappings on the Cantor fan.  We use our techniques to} obtain four different elementary proofs of Theorem \ref{main}:
\begin{theorem}\label{main}
{ 
There is a transitive homeomorphism on the Cantor fan.}
\end{theorem} 
Using closed relations on compact metric spaces, we also prove  the following theorem.
\begin{theorem}\label{mainn}{
There is a transitive function $f$ on a Cantor fan $X$ such that $\varprojlim(X,f)$ is a Lelek fan. In addition, the shift map on $\varprojlim(X,f)$ is a transitive homeomorphism.}
\end{theorem} 
{We believe that our techniques can be applied to produce other homeomorphisms and mappings to study dynamical properties of other compact metric spaces.  }

We proceed as follows. In Section \ref{s1}, the basic definitions and results that are needed later in the paper are presented. In Section \ref{s2}, we construct various dynamical systems using quotient spaces of topological products $X\times C$, where $X$ is a compact metric space and $C$ is a Cantor set. Here, we give our first two (of four) proofs of Theorem \ref{main}.  In Section \ref{s3},  we study dynamical systems with closed relations. 
Using this, we give two additional proofs of Theorem \ref{main}. In Section \ref{s4}, we study $\sigma$-transitivity of dynamical systems with closed relations. Then we prove Theorem \ref{mainn}.

\section{Definitions and Notation}\label{s1}
The following definitions, notation and well-known results are needed in the paper.

\begin{definition}
Let $X$ and $Y$ be metric spaces, and let $f:X\rightarrow Y$ be a function.  We use  $\Gamma(f)=\{(x,y)\in X\times Y \ | \ y=f(x)\}$
to denote \emph{  the graph of the function $f$}.
\end{definition}

\begin{definition}
Let $X$ be a metric space, $x\in X$ and $\varepsilon>0$. We use $B(x,\varepsilon)$ to denote the open ball,  {centered} at $x$ with radius $\varepsilon$.
\end{definition}
\begin{definition}
We use $\mathbb N$ to denote the set of positive integers and $\mathbb Z$ to denote the set of integers.  
\end{definition}
\begin{definition}
Let $X$ be a {non-empty} compact metric space and let ${F}\subseteq X\times X$ be a relation on $X$. If ${F}$ is closed in $X\times X$, then we say that ${F}$ is  \emph{  a closed relation on $X$}.  
\end{definition}

\begin{definition}
 \emph{A continuum} is a non-empty compact connected metric space.  \emph{A subcontinuum} is a subspace of a continuum, which is itself a continuum.
 \end{definition}
\begin{definition}
Let $X$ be a continuum.  We say that $X$ is \emph{a Cantor fan}, if $X$ is homeomorphic to the continuum $\bigcup_{c\in C}A_c$, where $C\subseteq [0,1]$ is a Cantor set and for each $c\in C$, $A_c$ is the  {convex} segment in the plane from $(0,0)$ to $(c,-1)$.
\begin{figure}[h!]
	\centering
		\includegraphics[width=15em]{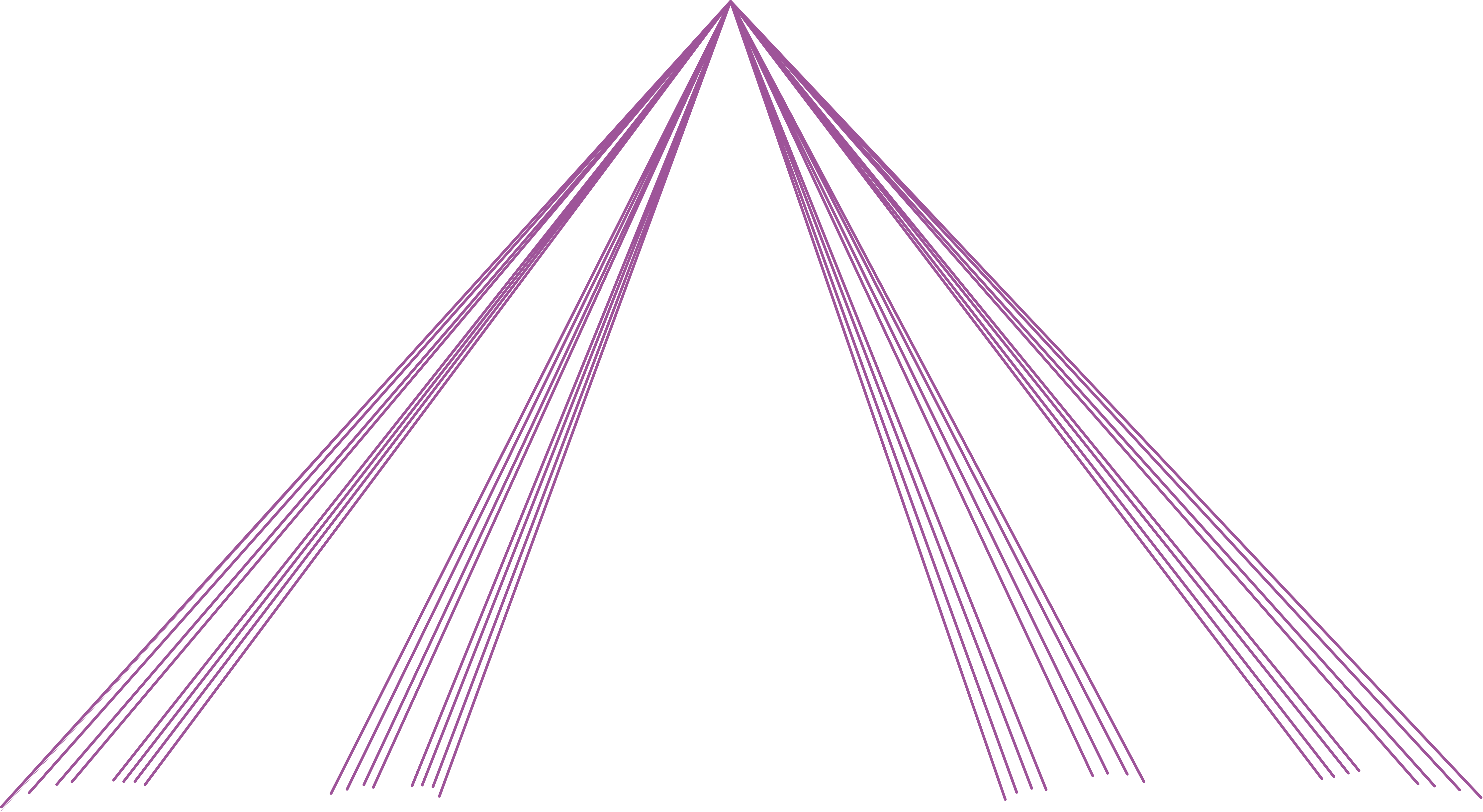}
	\caption{The Cantor fan}
	\label{fig000}
\end{figure}  
\end{definition}
\begin{definition}
Let $X$ be a Cantor fan and let $Y$ be a subcontinuum of $X$.  
A point $x\in Y$ is called an \emph{end-point of the continuum $Y$}, if for  every arc $A$ in $Y$ that contains $x$, $x$ is an end-point of $A$.  The set of all end-points of $Y$ will be denoted by $E(Y)$.
\end{definition}
\begin{definition} 
Let $X$ be a Cantor fan and let $Y$ be a subcontinuum of $X$.  We say that $Y$ is \emph{a Lelek fan}, if $\Cl(E(Y))=Y$.
\begin{figure}[h!]
	\centering
		\includegraphics[width=15em]{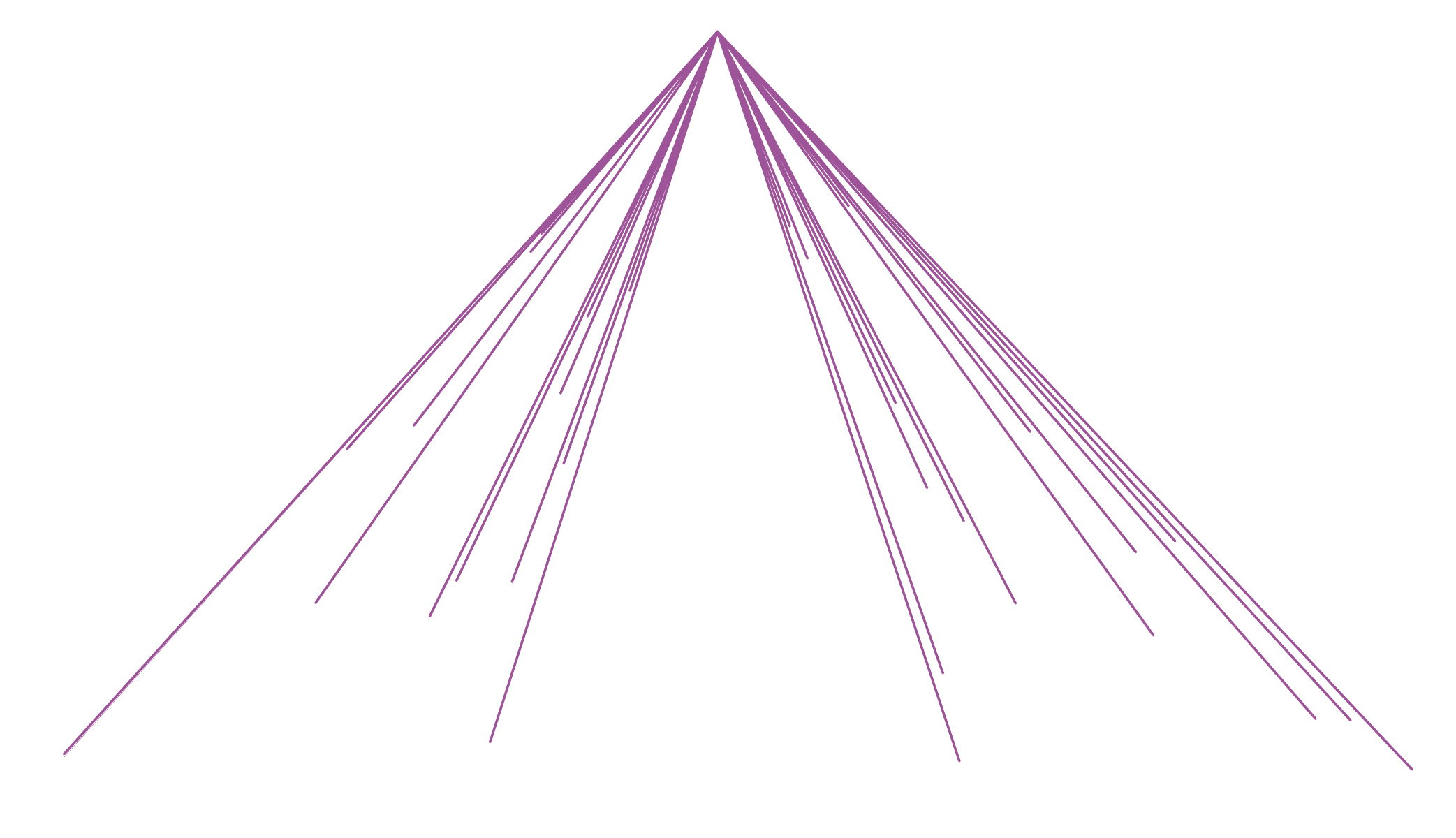}
	\caption{The Lelek fan}
	\label{figure2}
\end{figure}  
\end{definition}
The Lelek fan was constructed by A.~Lelek in \cite{lelek}.  An interesting  property of the Lelek fan $X$ is the fact that the set of its end-points is a dense one-dimensional set in $X$. It is also unique, i.e., it is the only non-degenerate smooth fan with a dense set of end-points.  This was proved independently by W.~D.~Bula and L.~Oversteegen  in \cite{oversteegen} and by W. ~Charatonik in  \cite{charatonik}. They proved that any two non-degenerate subcontinua of the Cantor fan with a dense set of endpoints are homeomorphic.  See \cite{nadler} for more information about continua, fans and their properties.
\begin{definition}
Let $X$ be a {non-empty} compact metric space and let ${F}$ be a closed relation on $X$. Then we call
$$
X_F^m=\Big\{(x_1,x_2,x_3,\ldots ,{x_m},x_{m+1})\in \prod_{{ i={1}}}^{{m{ +1}}}X \ | \ \textup{ for each } i\in{  \{{1,2},3,\ldots ,m\}}, (x_{i},x_{i+1})\in {F}\Big\}
$$
for each positive integer $m$, \emph{ the $m$-th Mahavier product of ${F}$}, { we call}
$$
X_F^+=\Big\{(x_1,x_2,x_3,\ldots )\in \prod_{{ i={1}}}^{\infty}X \ | \ \textup{ for each { positive} integer } i, (x_{i},x_{i+1})\in {F}\Big\}
$$
\emph{ the  Mahavier product of ${F}$}, and we call
$$
X_F=\Big\{(\ldots,x_{-3},x_{-2},x_{-1},{x_0}{ ;}x_1,x_2,x_3,\ldots )\in \prod_{i={-\infty}}^{\infty}X \ | \ \textup{ for each  integer } i, (x_{i},x_{i+1})\in {F}\Big\}
$$
\emph{ the two-sided  Mahavier product of ${F}$}.
\end{definition}

\begin{definition}
Let $X$ be a {non-empty} compact metric space and let ${F}$ be a closed relation on $X$. 
The function  $\sigma_F^{+} : {X_F^+} \rightarrow {X_F^+}$, 
 defined by 
$$
\sigma_F^{+} ({x_1,x_2,x_3,x_4},\ldots)=({x_2,x_3,x_4},\ldots)
$$
for each $({x_1,x_2,x_3,x_4},\ldots)\in {X_F^+}$, 
is called \emph{   the shift map on ${X_F^+}$}.      The function  $\sigma_F : {X_F} \rightarrow {X_F}$, 
 defined by 
$$
\sigma_F (\ldots,x_{-3},x_{-2},x_{-1},{x_0};x_1,x_2,x_3,\ldots )=(\ldots,x_{-3},x_{-2},x_{-1},{x_0},x_1;x_2,x_3,\ldots )
$$
for each $(\ldots,x_{-3},x_{-2},x_{-1},{x_0};x_1,x_2,x_3,\ldots )\in {X_F}$, 
is called \emph{   the shift map on ${X_F}$}.    
\end{definition}
\begin{observation}
Note that $\sigma_F$ is always a homeomorphism while $\sigma_F^+$ may not be a homeomorphism.
\end{observation}

\begin{definition}
Let $(X,f)$ be a dynamical system.  We say that $(X,f)$ is 
\begin{enumerate}
\item \emph{transitive}, if for all non-empty open sets $U$ and $V$ in $X$,  there is a non-negative integer $n$ such that $f^n(U)\cap V\neq \emptyset$.
\item \emph{dense orbit transitive},   if there is a point $x\in X$ such that its trajectory $\orbit_f(x)=\{x,f(x),f^2(x),f^3(x),\ldots\}$ is dense in $X$. We call such a point $x$ \emph{ a transitive point in $(X,f)$}.
\end{enumerate}
\end{definition}
\begin{definition}\label{povezava}
Let $(X,f)$ be a dynamical system. We say that the mapping $f$ is \emph{transitive}, if $(X,f)$ is transitive.
\end{definition}
\begin{observation}\label{isolatedpoints}
It is a well-known fact that if $X$ has no isolated points, then $(X,f)$ is transitive if and only if $(X,f)$ is dense orbit transitive. {Also, if $f$ is a homeomorphism, then $(X,f)$ is transitive if and only if $(X,f^{-1})$ is transitive.} See \cite{A,KS} for more information about transitive dynamical systems.
\end{observation}
\begin{definition}
Suppose $X$ is a compact metric space.  If $f:X \to X$ is a continuous function, the \emph{inverse limit space} generated by $f$ is the subspace
\begin{equation*}
 \varprojlim(X,f)=\Big\{(x_{1},x_{2},x_{3},\ldots ) \in \prod_{i=1}^{\infty} X \ | \ 
\text{ for each positive integer } i,x_{i}= f(x_{i+1})\Big\}
\end{equation*}
of the topological product $\prod_{i=1}^{\infty} X$.  { The function  $\sigma : \varprojlim(X,f) \rightarrow \varprojlim(X,f)$, 
 defined by 
$$
\sigma (x_1,x_2,x_3,x_3,\ldots )=(x_2,x_3,x_4,\ldots )
$$
for each $(x_1,x_2,x_3,\ldots )\in \varprojlim(X,f)$, 
is called \emph{   the shift map on $\varprojlim(X,f)$}.    }
\end{definition}
The following theorem is a well-known result proved in  \cite{he},  and it also easily follows  from the results from \cite{li}.  
\begin{theorem}\label{shifttransitive}
Let $(X,f)$ be a dynamical system. If  $(X,f)$ is transitive and $\sigma :\varprojlim(X,f)\rightarrow \varprojlim(X,f)$ is the shift map on $\varprojlim(X,f)$, then also $(\varprojlim(X,f),\sigma)$ is a transitive dynamical system.  
\end{theorem}
\begin{observation}\label{homeomorphism}
For any dynamical system $(X,f)$, the shift map  $\sigma$ on $\varprojlim(X,f)$ is a homeomorphism from $\varprojlim(X,f)$ to $\varprojlim(X,f)$.
\end{observation}

\section{Transitive homeomorphisms on quotient spaces}\label{s2}
In this section, we give our first two proofs of Theorem \ref{main}.  We start the section with some well-known results that are needed later. To make the paper more reader friendly, we give their proofs (since the proofs are elementary and short).
\begin{definition}
Let $X$ be a compact metric space and let $\sim$ be an equivalence relation on $X$. We use 
\begin{enumerate}
\item $[x]$ to denote the equivalence class $\{y\in X \ | \ y\sim x\}$ of the element $x$ for any $x\in X$,
\item $X/_{\sim}$ to denote the quotient space $\{[x] \ | \ x\in X\}$, and
\item $q$ to denote the natural quotient mapping $X\rightarrow X/_{\sim}$, which is defined by $q(x)=[x]$ for any $x\in X$.   
\end{enumerate}
\end{definition}
\begin{observation}\label{opiopi}
Let $X$ be a compact metric space and let $\sim$ be an equivalence relation on $X$. Then the following hold.
\begin{enumerate}
\item The quotient map $q$ is a continuous surjection.
\item For any $U\subseteq X/_{\sim}$, 
$$
U \textup{ is open in } X/_{\sim}  ~~~  \Longleftrightarrow  ~~~  q^{-1}(U) \textup{ is open in } X.
$$
\item If the quotient space  $X/_{\sim}$ is a Hausdorff space, then it is a metrizable compactum. Therefore, in this case,  $X/_{\sim}$ is a continuum, if it is connected.
\end{enumerate}
\end{observation}
\begin{definition}
Let $X$ be a compact metric space, let $\sim$ be an equivalence relation on $X$,  and let $f:X\rightarrow X$ be a function such that for all $x,y\in X$,
$$
x\sim y  \Longleftrightarrow f(x)\sim f(y).
$$
 Then we let $f^{\star}:X/_{\sim}\rightarrow X/_{\sim}$ be defined by   
$
f^{\star}([x])=[f(x)]
$
for any $x\in X$. 
\end{definition}
The following theorem is a well-known result. Since its proof is short, we give it for the completeness of the paper.
\begin{theorem}\label{kvocienti}
Let $X$ be a compact metric space, let $\sim$ be an equivalence relation on $X$, and  let $f:X\rightarrow X$ be a function such that for all $x,y\in X$,
$$
x\sim y  \Longleftrightarrow f(x)\sim f(y).
$$
 Then the following hold.
\begin{enumerate}
\item $f^{\star}$ is a well-defined function from  $X/_{\sim}$ to $X/_{\sim}$. 
\item If $f$ is continuous, then $f^{\star}$ is continuous.
\item If $f$ is a homeomorphism, then $f^{\star}$ is a homeomorphism.
\item If $f$ is transitive, then $f^{\star}$ is transitive.
\end{enumerate}
\end{theorem}
\begin{proof}
To see that $f^{\star}$ is a well-defined function from  $X/_{\sim}$ to $X/_{\sim}$, let $x,y\in X$ be such that $x\sim y$.  It follows that
$$
f^{\star}([x])=[f(x)]=[f(y)]=f^{\star}([y]).
$$
Next, suppose that $f$ is continuous. To show that  $f^{\star}$ is continuous, let $U$ be an open set in $X/_{\sim}$. We show that $(f^{\star})^{-1}(U)$ is open in $X/_{\sim}$ by showing that $q^{-1}((f^{\star})^{-1}(U))$ is open in $X$. Since $f^{\star}\circ q=q\circ f$, it follows that 
$$
q^{-1}((f^{\star})^{-1}(U))=(f^{\star}\circ q)^{-1}(U)=(q\circ f)^{-1}(U).
$$
Since $f$ and $q$ are continuous, it follows that $q^{-1}((f^{\star})^{-1}(U))$ is open in $X$.  Therefore,  $f^{\star}$ is continuous. Now, suppose that $f$ is a homeomorphism.  To see that $f^{\star}$ is bijective,  let $g:X/_{\sim}\rightarrow X/_{\sim}$ be defined by 
$
g([x])=[f^{-1}(x)]
$
for any $x\in X$.
It follows that for each $x\in X$,
\begin{align*}
f^{\star}(g([x]))=&f^{\star}([f^{-1}(x)])=[f(f^{-1}(x))]=[x],
\\
g(f^{\star}([x]))=&g([f(x)])=[f^{-1}(f(x))]=[x].
\end{align*}
Therefore, $f^{\star}$ is bijective and $(f^{\star})^{-1}=g$.  The proof  that $g$ is continuous is similar to the proof that $f^{\star}$ is continuous. We leave the details to the reader.  Since $f^{\star}$ and $(f^{\star})^{-1}$ are both continuous, it follows that $f^{\star}$ is a homeomorphism.
Finally, suppose that $f$ is transitive.  Let $U$ and $V$ be any non-empty open sets in $X/_{\sim}$.  We prove that there is a non-negative integer $n$ such that $(f^{\star})^n(U)\cap V\neq \emptyset$.  Since the quotient map is a continuous surjection, it follows that $q^{-1}(U)$ and $q^{-1}(V)$ are non-empty open subsets of $X$.  Let $n$ be a non-negative integer such that $f^n(q^{-1}(U))\cap q^{-1}(V)\neq \emptyset$ and let $x\in f^n(q^{-1}(U))\cap q^{-1}(V)$.  We show that $[x]\in (f^{\star})^n(U)\cap V$. It follows from $x\in q^{-1}(V)$ that $q(x)\in q(q^{-1}(V))$ and, therefore, $[x]\in V$.  Also, it follows from $x\in f^n(q^{-1}(U))$ that $q(x)\in q(f^n(q^{-1}(U)))$ and, therefore, $[x]\in q(f^n(q^{-1}(U)))$. To show that also $[x]\in (f^{\star})^n(U)$, we prove that $q(f^n(q^{-1}(U)))=(f^{\star})^n(U)$:
\begin{align*}
(f^{\star})^n(U)=&\{(f^{\star})^n([t]) \ | \ [t]\in U\}=\{[f^n(t)] \ | \ q(t)\in U\}=
\\
=&\{q(f^n(t)) \ | \ q(t)\in U\}=q(f^n(\{t \ | \ q(t)\in U\}))=q(f^n(q^{-1}(U))).
\end{align*}
This completes the proof.
\end{proof}
Let $X$ be a set and let $n$ be a positive integer.  For each positive integer $i$ and for each point $\mathbf x=(x_1,x_2,x_3,\ldots)\in \prod_{k=1}^{\infty}X$ (or for each $i\in \{1,2,3,\ldots, n\}$ and for each point $\mathbf x=(x_1,x_2,x_3,\ldots,x_n)\in \prod_{k=1}^{n}X$ or for each integer $i$ and for each point $\mathbf x=({\ldots,x_{-2},x_{-1}},{x_0;x_1,x_2,\ldots})\in \prod_{k=-\infty}^{\infty}X$), we use  $\mathbf x(i)$  to denote the $i$-th coordinate $x_i$ of the point $\mathbf x$.
For all positive  integers $i$ and $j$ (or for all $i,j\in \{1,2,3,\ldots, n\}$ or for all integers $i$ and $j$) such that $i\leq j$,  we use $[i,j]$ to denote the set $\{i, i+1,i+2,\ldots,j\}$ and $\pi_{[i,j]}$ to denote the projection $\pi_{[i,j]}:\prod_{k=1}^{\infty}X\rightarrow \prod_{k=i}^{j}X$ (or the projection $\pi_{[i,j]}:\prod_{k=1}^{n}X\rightarrow \prod_{k=i}^{j}X$ or the projection $\pi_{[i,j]}:\prod_{k=-\infty}^{\infty}X\rightarrow \prod_{k=i}^{j}X$) that is defined by
$
\pi_{[i,j]}(x_1,x_2,x_3,\ldots)=(x_i,x_{i+1},x_{i+2},\ldots, x_{j})
$
for any $(x_1,x_2,x_3,\ldots)\in \prod_{k=1}^{\infty}X$ (or 
$
\pi_{[i,j]}(x_1,x_2,x_3,\ldots,x_n)=(x_i,x_{i+1},x_{i+2},\ldots, x_{j})
$
for any $(x_1,x_2,x_3,\ldots,x_n)\in \prod_{k=1}^{n}X$ or 
$
\pi_{[i,j]}({\ldots,x_{-1},x_0};{x_1,x_2,\ldots})=(x_i,x_{i+1},\ldots, x_{j})
$
for any $({\ldots,x_{-1},x_0};{x_1,x_2,\ldots})\in \prod_{k=-\infty}^{\infty}X$).
For $i=j$, we use $\pi_i$ to denote the projection  $\pi_{[i,i]}$.

For a {non-empty} compact metric space $X$,   we use 
$p_1:X\times X\rightarrow X$ and $p_2:X\times X\rightarrow X$
to denote \emph{ the standard projections} defined by
$p_1(s,t)=s$ and $p_2(s,t)=t$ for all $(s,t)\in X\times X$.
\begin{definition}
For each positive integer $n$, we use 
\begin{enumerate}
\item $C_n$ to denote the topological product 
$
C_n=\prod_{k=1}^{\infty}\{1,2,3,\ldots ,n\},
$
 where the set $\{1,2,3,\ldots, n\}$ is equipped with the discrete topology.  We also use $\sigma_n$ to denote the shift map $\sigma_n:C_n\rightarrow C_n$, defined by 
$$
{ \sigma_n(\mathbf x)}=\sigma_n(\mathbf x(1),\mathbf x(2),\mathbf x(3),\ldots)=(\mathbf x(2),\mathbf x(3),\mathbf x(4),\ldots)
$$
for any $\mathbf x\in C_n$.
\item $D_n$ to denote the topological product 
$
D_n=\prod_{k=-\infty}^{\infty}\{1,2,3,\ldots ,n\},
$
 where the set $\{1,2,3,\ldots, n\}$ is equipped with the discrete topology.  We also use $\tau_n$ to denote the shift map $\tau_n:D_n\rightarrow D_n$, defined by 
$$
{ \tau_n(\mathbf x)}=\tau_n({\ldots,\mathbf x(-1),\mathbf x(0)};{\mathbf x(1),\mathbf x(2),\ldots})=({\ldots,\mathbf x(-1),\mathbf x(0),\mathbf x(1)};{\mathbf x(2)},\ldots)
$$
for any $\mathbf x\in D_n$.
\end{enumerate}
\end{definition}
\begin{observation}
Note that for each positive integer $n\geq 2$, $C_n$ and $D_n$ {(equipped with the product topologies)} are Cantor sets and that the shift maps $\sigma_n$ and $\tau_n$ are transitive. Also, note that for each positive integer $n$,  $\tau_n$ is a homeomorphism while $\sigma_n$ is not. 
\end{observation}
\begin{definition}
Let $X$ be a compact metric space, let $n$ be a positive integer, and let $f_1,f_2,f_3,\ldots,f_n:X\rightarrow X$ be any functions.  Then we use
\begin{enumerate}
\item $(f_1,\ldots,f_n)_{C_n}$ to denote the mapping from $C_n\times X$ to $C_n\times X$, which is defined by 
$$
(f_1,\ldots,f_n)_{C_n}(\mathbf x,t)=(\sigma_n(\mathbf x), f_{\mathbf x(1)}(t))
$$
for any $(\mathbf x,t)\in C_n\times X$.  {  We call $(f_1,\ldots,f_n)_{C_n}$ \emph{a sorting function  with respect to the first projection of $C_n$}.}
\item $(f_1,\ldots,f_n)_{D_n}$ to denote the mapping from $D_n\times X$ to $D_n\times X$, which is defined by 
$$
(f_1,\ldots,f_n)_{D_n}(\mathbf x,t)=(\tau_n(\mathbf x), f_{\mathbf x(1)}(t))
$$
for any $(\mathbf x,t)\in D_n\times X$. {  We call $(f_1,\ldots,f_n)_{D_n}$ \emph{a sorting function  with respect to the first projection of  $D_n$}.}
\end{enumerate} 
\end{definition}
{  In Example \ref{sorting},  we demonstrate how the defined sorting functions work. 
\begin{example}\label{sorting}
Let $n=2$, let $X=[0,1]$,  and let $f_1,f_2:[0,1]\rightarrow [0,1]$ be defined by $f_1(t)=t$ and $f_2(t)=t^2$ for any $t\in [0,1]$. Then for $\mathbf x=(1,2,1,2,1,\ldots)$,  $\mathbf x(1)=1$ and it follows that for each $t\in [0,1]$,
$$
(f_1,f_2)_{C_2}(\mathbf x,t)=(\sigma_2(\mathbf x), f_{\mathbf x(1)}(t))=(\sigma_2(1,2,1,2,1,\ldots), f_1(t))=((2,1,2,1,\ldots), t^2),
$$
and for $\mathbf x=(\ldots, 2,1;2,1,\ldots)$,  $\mathbf x(1)=2$ and it follows that for each $t\in [0,1]$,
$$
(f_1,f_2)_{D_2}(\mathbf x,t)=(\tau_2(\mathbf x), f_{\mathbf x(1)}(t))=(\tau_2(\ldots, 2,1;2,1,\ldots)), f_2(t))=((\ldots, 2,1,2;1,\ldots)), t).
$$
\end{example}
}
\begin{proposition}\label{propro}
Let $X$ be a compact metric space, let $n$ be a positive integer, and let $f_1,f_2,f_3,\ldots,f_n:X\rightarrow X$ be any functions.  Then the following hold.
\begin{enumerate}
\item \label{1} If $f_1$, $f_2$, $f_3$,  $\ldots$, $f_{n}$ are continuous, then $(f_1,\ldots,f_n)_{C_n}$ and $(f_1,\ldots,f_n)_{D_n}$ are continuous.
\item \label{2} If $f_1$, $f_2$, $f_3$,  $\ldots$, $f_{n}$  are surjective,  then $(f_1,\ldots,f_n)_{C_n}$ and $(f_1,\ldots,f_n)_{D_n}$ are surjective.
\item If $f_1$, $f_2$, $f_3$,  $\ldots$,  $f_n$ are homeomorphisms,  then $(f_1,\ldots,f_n)_{D_n}$ is a homeomorphism.
\end{enumerate}
\end{proposition}
\begin{proof}
Suppose that the functions $f_1$, $f_2$, $f_3$,  $\ldots$, $f_n$ are continuous. First, we prove that $(f_1,\ldots,f_n)_{C_n}$ is continuous.  {Let $\delta=1$ and let}
 $U=V\times W$, where  $W$ is an open set in $X$ and $V$ is an open set in $C_n$ such that $\diam(V)<\delta$.  Also there is $k\in \{1,2,3,\ldots,n\}$ such that $\pi_1(V)=\{k\}$. Fix such an element $k$. To see that $(f_1,\ldots,f_n)_{C_n}$ is continuous, we prove that $(f_1,\ldots,f_n)_{C_n}^{-1}(U)$ is open in $C_n\times X$.  Observe that 
$$
(f_1,\ldots,f_n)_{C_n}^{-1}(U)=(f_1,\ldots,f_n)_{C_n}^{-1}(V\times W)=\sigma_n^{-1}(V)\times f_{k}^{-1}(W).
$$
Since the functions $\sigma_n$ and $f_k$ are continuous, it follows that $(f_1,\ldots,f_n)_{C_n}^{-1}(U)$ is open in $C_n\times X$.  

{  The proof that $(f_1,\ldots,f_n)_{D_n}$ is continuous is analogous to the proof above - we leave it for the reader. }

Now, suppose that $f_1$, $f_2$, $f_3$,  $\ldots$,  $f_n$ are surjective.  To show that $(f_1,\ldots,f_n)_{C_n}$ is surjective, let $(\mathbf x,t)\in C_n\times X$ be any point. Since $\sigma_n$ is surjective, there is a point $\mathbf t \in C_n$ such that $\sigma_n(\mathbf t)=\mathbf x$. Choose such a point $\mathbf t$. Then
$$
(f_1,\ldots,f_n)_{C_n}(\mathbf t,f_{\mathbf t(1)}^{-1}(t))=(\sigma_n(\mathbf t), f_{\mathbf t(1)}(f_{\mathbf t(1)}^{-1}(t)))=(\mathbf x,t),
$$
therefore, $(f_1,\ldots,f_n)_{C_n}$ is surjective. 

{  The proof that $(f_1,\ldots,f_n)_{D_n}$ is surjective is analogous to the proof above that $(f_1,\ldots,f_n)_{C_n}$ is surjective - we leave it for the reader. }

Finally, suppose that $f_1$, $f_2$, $f_3$,  $\ldots$,  $f_n$ are homeomorphisms.  It follows from \ref{1}  that $(f_1,\ldots,f_n)_{D_n}$ is a continuous function.  To see that $(f_1,\ldots,f_n)_{D_n}$ is a bijection, let $h=(f_1,\ldots,f_n)_{D_n}$ and let $g:D_n\times X\rightarrow D_n\times X$ be defined by 
$$
g(\mathbf x,t)=(\tau_n^{-1}(\mathbf x), f_{\tau_n^{-1}(\mathbf x)(1)}^{-1}(t))
$$
for any $(\mathbf x,t)\in D_n\times X$. Observe that for each $(\mathbf x,t)\in D_n\times X$,
\begin{align*}
g(h(\mathbf x,t))&=g(\tau_n(\mathbf x), f_{\mathbf x(1)}(t))=(\tau_n^{-1}(\tau_n(\mathbf x)), f_{\tau_n^{-1}(\tau_n(\mathbf x))(1)}^{-1}(f_{\mathbf x(1)}(t)))=(\mathbf x,t),
\\
h(g(\mathbf x,t))&=h(\tau_n^{-1}(\mathbf x), f_{\tau_n^{-1}(\mathbf x)(1)}^{-1}(t))=(\tau_n(\tau_n^{-1}(\mathbf x)), f_{\tau_n^{-1}(\mathbf x)(1)}(f_{\tau_n^{-1}(\mathbf x)(1)}^{-1}(t)))=(\mathbf x,t).
\end{align*}
It follows that $(f_1,\ldots,f_n)_{D_n}$ is a bijection. 
Since $(f_1,\ldots,f_n)_{D_n}$ is a continuous bijection from a compact metric space onto itself, it is a homeomorphism.
\end{proof}
In the following subsections, we examine separately the quotient spaces $C_n\times X$ and $D_n\times X$. We will use the notion of different types of concatenation that is defined in the following definition.

\begin{definition}
Let $X$ be a compact metric space, let $n$ and $m$ be positive integers, and let $\mathbf x\in \prod_{k=1}^{n}X$,  $\mathbf y\in \prod_{k=1}^{m}X$, $\mathbf z\in \prod_{i=1}^{\infty}X$ and $\mathbf w\in \prod_{i=-\infty}^{0}X$. Then we define 
\begin{enumerate}
\item \emph{the concatenation $\mathbf x\oplus \mathbf y\in \prod_{k=1}^{m+n}X$} of $\mathbf x$ and $\mathbf y$ as follows: 
$$
\mathbf x\oplus \mathbf y=(\mathbf x(1),\mathbf x(2),\mathbf x(3),\ldots,\mathbf x(n),\mathbf y(1),\mathbf y(2),\mathbf y(3),\ldots,\mathbf y(m)).
$$
\item  \emph{the concatenation $\mathbf x\oplus \mathbf z\in \prod_{k=1}^{\infty}X$} of $\mathbf x$ and $\mathbf z$ as follows: 
$$
\mathbf x\oplus \mathbf z=(\mathbf x(1),\mathbf x(2),\mathbf x(3),\ldots,\mathbf x(n),\mathbf z(1),\mathbf z(2),\mathbf z(3),\ldots).
$$
\item  \emph{the concatenation $\mathbf w\oplus \mathbf x\in \prod_{k=-\infty}^{0}X$} of $\mathbf w$ and $\mathbf x$ as follows: 
$$
{\mathbf w}\oplus {\mathbf x}=({\ldots, \mathbf w(-3),\mathbf w(-2),\mathbf w(-1), \mathbf w(0)},{\mathbf x(1),\mathbf x(2),\mathbf x(3),\ldots, \mathbf x(n)}).
$$
\item  \emph{the concatenation $\mathbf w\oplus \mathbf z\in \prod_{k=-\infty}^{\infty}X$} of $\mathbf w$ and $\mathbf z$ as follows: 
$$
{\mathbf w}\oplus {\mathbf z}=({\ldots, \mathbf w(-3),\mathbf w(-2),\mathbf w(-1),\mathbf w(0)};{\mathbf z(1),\mathbf z(2),\mathbf z(3),\ldots}).
$$
\end{enumerate}
\end{definition}
\subsection{The quotient spaces of  $C_n\times X$ }
\begin{definition}
Let $X$ be a compact metric space, let $n$ be a positive integer, let 
$\mathbf x\in \prod_{k=1}^{n}X$ and $\mathbf y\in \prod_{k=1}^{\infty}X$, and let $\varepsilon >0$. We say that 
$
\big|\mathbf x-\mathbf y\big|_n<\varepsilon,
$
if
\begin{enumerate}
\item for each $k\in \{1,2,3,\ldots,n\}$, $\mathbf x(k)=\mathbf y(k)$ and
\item $\displaystyle \sum_{k=n}^{\infty}\frac{\diam(X)}{2^k}<\varepsilon$.
\end{enumerate}
\end{definition}

\begin{definition}
Let $X$ be a compact metric space, let $n$ be a positive integer, and let $f_1,f_2,f_3,\ldots,f_n:X\rightarrow X$ be any functions.  Also,  let $i$ and $j$ be  integers such that $i\leq j$  and let ${\mathbf x}\in \prod_{k=i}^{j}\{1,2,3,\ldots,n\}$.  Then we use $f_{{\mathbf x}}$ to denote the function $f_{{\mathbf x}}:X\rightarrow [0,1]$,  defined by
$$
f_{{\mathbf x}}(t)=(f_{{\mathbf x}(j)}\circ f_{{\mathbf x}(j-1)}\circ f_{{\mathbf x}(j-2)}\circ \ldots \circ f_{{\mathbf x}(i+2)}\circ f_{{\mathbf x}(i+1)}\circ f_{{\mathbf x}(i)})(t)
$$
for any $t\in X$.
\end{definition}
\begin{definition}
Let $(X,d)$ be a compact metric space, let $n$ be a positive integer, and let $f_1,f_2,f_3,\ldots,f_n:X\rightarrow X$ be any functions.  Also, let $U$ be a dense open set in $X$.  We say that  \emph{ the $n$-tuple $(f_1,f_2, f_3, \ldots,f_n)$ has property$\mathcal L(U)$}, if for each $\varepsilon >0$, for each $x\in U$ and for each $z\in X$, there are  non-negative integers $k_1$, $k_2$, $k_3$, $\ldots$,  $k_n$ such that
$$
d((f_1^{k_1}\circ f_2^{k_2}\circ f_3^{k_3} \circ \ldots \circ f_n^{k_n})(x),z)<\varepsilon.
$$
We say that  \emph{the $n$-tuple $(f_1,f_2, f_3, \ldots,f_n)$ has property$\mathcal L$},  if there is a dense open set $U$ in $X$ such that  the $n$-tuple $(f_1,f_2, f_3, \ldots,f_n)$ has property$\mathcal L(U)$.
\end{definition}

\begin{theorem}\label{glav}
Let $X$ be a compact metric space that does not have any isolated points, let $n$ be a positive integer, and let $f_1,f_2,f_3,\ldots,f_n:X\rightarrow X$ be any homeomorphisms.  If the $n$-tuple $(f_1,f_2, f_3, \ldots,f_n)$ has property $\mathcal L$,  then  $(f_1,\ldots,f_n)_{C_n}$ is a transitive continuous surjection.
\end{theorem}
\begin{proof}
Let $U$ be a dense open set in $X$ such that the $n$-tuple $(f_1,f_2, f_3, \ldots,f_n)$ has property $\mathcal L(U)$.
To  prove that $(f_1,\ldots,f_n)_{C_n}$ is a transitive continuous surjection, let  $h=(f_1,\ldots,f_n)_{C_n}$.  It follows from Proposition \ref{propro} that $h$ is a continuous surjection.  We show that $h$ is transitive by constructing a transitive point in $(C_n\times X,h)$. 
Let $\{(\mathbf x_i,t_i) \ | \ i \textup{ is a positive integer}\}$ be a countable dense subset of $C_n\times X$. Also,  for each positive integer $i$,  let $\varepsilon_i>0$ be such that $\displaystyle \lim_{i\to \infty}\varepsilon_i=0$. For each positive integer $i$, let $m_i$ be a positive integer and let $\hat{\mathbf x}_i\in \prod_{k=1}^{m_i}\{1,2,3,\ldots,n\}$ be such that
$
\big| \hat{\mathbf x}_i-\mathbf x_i\big|_{m_i}<\varepsilon_i.
$
For each $x\in U$, for each $z\in X$ and for each $\varepsilon >0$, let $k_1(x,z,\varepsilon)$, $k_2(x,z,\varepsilon)$, $k_3(x,z,\varepsilon)$, $\ldots$,  $k_n(x,z,\varepsilon)$ be positive integers such that 
$$
d((f_1^{k_1(x,z,\varepsilon)}\circ f_2^{k_2(x,z,\varepsilon)}\circ f_3^{k_3(x,z,\varepsilon)} \circ \ldots \circ f_n^{k_n(x,z,\varepsilon)})(x),z)<\varepsilon,
$$
and let 
$$
Y(x,z,\varepsilon)=(\underbrace{n,n,\ldots ,n}_{k_n(x,z,\varepsilon)},\underbrace{n-1,n-1,\ldots ,n-1}_{k_{n-1}(x,z,\varepsilon)}, \underbrace{n-2,n-2,\ldots ,n-2}_{k_{n-2}(x,z,\varepsilon)},\ldots , \underbrace{1,1,\ldots ,1}_{k_1(x,z,\varepsilon)}).
$$
Now, we inductively define points ${\mathbf w}_i$ (based on $i$ being even or odd) as follows. Let $x_0\in U$ and let 
$
\hat{\mathbf w}_1=Y(x_0,t_1,\varepsilon_1),
$
let 
$
\hat{\mathbf w}_2=\hat{\mathbf w}_1\oplus \hat{\mathbf x}_1
$
and let
$
\hat{\mathbf w}_3=\hat{\mathbf w}_2\oplus Y\Big(f_{\hat{\mathbf w}_2}(x_0),t_2,\varepsilon_2\Big).
$
Let $i$ be any positive integer and suppose that $\hat{\mathbf w}_1$, $\hat{\mathbf w}_2$, $\hat{\mathbf w}_3$, $\ldots$, $\hat{\mathbf w}_{2i-1}$ have already been defined.  Then we define $\hat{\mathbf w}_{2i}$ and $\hat{\mathbf w}_{2i+1}$ as follows:
$
\hat{\mathbf w}_{2i}=\hat{\mathbf w}_{2i-1}\oplus \hat{\mathbf x}_i
$
and
$
\hat{\mathbf w}_{2i+1}=\hat{\mathbf w}_{2i}\oplus Y\Big(f_{\hat{\mathbf w}_{2i}}(x_0),t_{i+1},\varepsilon_{i+1}\Big).
$
Next, for each positive integer $i$, let 
$
\mathbf w_i=\hat{\mathbf w}_{i}\oplus (1,1,1,\ldots).
$
Note that the sequence $(\mathbf w_i)$ is convergent in $C_n$ and let
$
\mathbf w=\lim_{i\to \infty}\mathbf w_i.
$
Finally, we show that $(\mathbf w,x_0)$ is a transitive point in $(C_n\times X,h)$.
Let $V$ be a non-empty open set in $C_n$, let $W$ be a non-empty open set in $X$, and let $U_0=V\times W$.  We show that there is a positive integer $k$ such that $h^k(\mathbf w,x_0)\in U_0$.  Let $k_0$ be a positive integer such that $B((\mathbf x_{k_0},t_{k_0}),\varepsilon_{k_0})\subseteq U_0$.   Also, let $i$ and $k$ be the positive integers such that 
$
\pi_{[1,i]}(\sigma_n^{k}(\mathbf w))=\hat{\mathbf x}_{k_0}.
$
Then $h^k(\mathbf w,x_0)\in U_0$ and we are done.
\end{proof}
{ 
\begin{observation}
Note that $(f_1,\ldots,f_n)_{C_n}$ from Theorem \ref{glav} is not a homeomorphism. This  follows from the fact that $\sigma_n$ is not a homeomorphism.  
\end{observation}}
\begin{definition}
Let $m$ and $n$ be positive integers,  and let $\mathbf y\in C_n$.   We use $\mathbf y[m]$ to denote 
$
\mathbf y[m]=(\mathbf y(m),\mathbf y(m-1),\mathbf y(m-2),\ldots, \mathbf y(1)).
$
\end{definition}
\begin{definition}
Let $X$ be a compact metric space, let $m$ and $n$ be positive integers,   let $\mathbf y\in C_n$,  and let $f_1,f_2,f_3,\ldots,f_n:X\rightarrow X$ be any homeomorphisms.   We use $f_{\mathbf y[m]}^{-1}$ to denote the function $f_{\mathbf y[m]}^{-1}:[0,1]\rightarrow [0,1]$ defined by 
$$
f_{\mathbf y[m]}^{-1}(t)=\Big(f_{\mathbf y(m)}^{-1}\circ f_{\mathbf y(m-1)}^{-1}\circ f_{\mathbf y(m-2)}^{-1}\circ \ldots \circ f_{\mathbf y(1)}^{-1}\Big)(t)
$$
for each $t\in [0,1]$.
\end{definition}
\begin{observation}
Let $X$ be a compact metric space, let  $n$ be a positive integer,  and let $f_1,f_2,f_3,\ldots,f_n:X\rightarrow X$ be any homeomorphisms.  
Note that for each $\mathbf x\in \varprojlim(C_n\times X,(f_1,\ldots,f_n)_{C_n})$, there are uniquely determined points $\mathbf a,\mathbf b\in C_n$ and $t\in X$ such that 
$$
\mathbf x=\Big(\big(\mathbf a,t\big),\big(\mathbf b[1]\oplus \mathbf a,f_{\mathbf b[1]}^{-1}(t)\big),\big(\mathbf b[2]\oplus \mathbf a,f_{\mathbf b[2]}^{-1}(t)\big),\big(\mathbf b[3]\oplus \mathbf a,f_{\mathbf b[3]}^{-1}(t)\big),\ldots\Big).
$$
\end{observation}
\begin{theorem}\label{loj}
Let $X$ be a compact metric space, let  $n$ be a positive integer,  and let $f_1,f_2,f_3,\ldots,f_n:X\rightarrow X$ be any homeomorphisms.  
Then the  inverse limit space $\varprojlim(C_n\times X,(f_1,\ldots,f_n)_{C_n})$ is homeomorphic to $C_n\times X$.
\end{theorem}
\begin{proof}
Note that $C_n\times X$ is homeomorphic to $(C_n\times C_n)\times X$. We show that  the inverse limit $\varprojlim(C_n\times X,(f_1,\ldots,f_n)_{C_n})$ is homeomorphic to $(C_n\times C_n)\times X$.  For each $\mathbf x\in \varprojlim(C_n\times X,(f_1,\ldots,f_n)_{C_n})$, let $\mathbf a_{\mathbf x},\mathbf b_{\mathbf x}\in C_n$ and $t_{\mathbf x}\in X$ be such that
$$
\mathbf x=\Big(\big(\mathbf a_{\mathbf x},t_{\mathbf x}\big),\big(\mathbf b_{\mathbf x}[1]\oplus \mathbf a_{\mathbf x},f_{\mathbf b_{\mathbf x}[1]}^{-1}(t_{\mathbf x})\big),\big(\mathbf b_{\mathbf x}[2]\oplus \mathbf a_{\mathbf x},f_{\mathbf b_{\mathbf x}[2]}^{-1}(t_{\mathbf x})\big),\big(\mathbf b_{\mathbf x}[3]\oplus \mathbf a_{\mathbf x},f_{\mathbf b_{\mathbf x}[3]}^{-1}(t_{\mathbf x})\big),\ldots\Big).
$$
Let $T:\varprojlim(C_n\times X,(f_1,\ldots,f_n)_{C_n})\rightarrow (C_n\times C_n)\times X$ be defined by 
$
T(\mathbf x)=((\mathbf a_{\mathbf x},\mathbf b_{\mathbf x}),t_{\mathbf x})
$
for each $\mathbf x\in \varprojlim(C_n\times X,(f_1,\ldots,f_n)_{C_n})$. Clearly, $T$ is a homeomorphism.
\end{proof}
\begin{definition}
Let $X$ be a compact metric space, let  $n$ be a positive integer,  and let $f_1,f_2,f_3,\ldots,f_n:X\rightarrow X$ be any homeomorphisms.  
We use $\sigma_{(f_1,\ldots,f_n)_{C_n}}$ to denote the shift homeomorphism on $\varprojlim(C_n\times X,(f_1,\ldots,f_n)_{C_n})$.
\end{definition}
\begin{theorem}\label{miladojka1}
Let $X$ be a compact metric space, let  $n$ be a positive integer,  and let $f_1,f_2,f_3,\ldots,f_n:X\rightarrow X$ be any homeomorphisms.  If  the $n$-tuple $(f_1, f_2, f_3,\ldots, f_n)$ has property$\mathcal L$, then the shift homeomorphism $\sigma_{(f_1,\ldots,f_n)_{C_n}}$ is transitive.
\end{theorem}
\begin{proof}
By Theorem \ref{glav},  the function $(f_1,\ldots,f_n)_{C_n}$ is transitive.  It follows from Theorem \ref{shifttransitive} that $\sigma_{(f_1,\ldots,f_n)_{C_n}}$ is transitive. 
\end{proof}
\begin{definition}
Let $X$ be a compact metric space, let  $n$ be a positive integer,  let $f_1,f_2,f_3,\ldots,f_n:X\rightarrow X$ be any homeomorphisms, and let  $\sim$ be an equivalence relation on $(C_n\times C_n)\times X$.  Then we define a relation $\sim_{(f_1,\ldots,f_n)_{C_n}}$ on the inverse limit $\varprojlim(C_n\times X,(f_1,\ldots,f_n)_{C_n})$ by
$$
\mathbf x\sim_{(f_1,\ldots,f_n)_{C_n}} \mathbf y ~~~  \Longleftrightarrow ~~~ \mathbf x=\mathbf y  \textup{ or }  ((\mathbf a_{\mathbf x},\mathbf b_{\mathbf x}),t_{\mathbf x})\sim ((\mathbf a_{\mathbf y},\mathbf b_{\mathbf y}),t_{\mathbf y})
$$
for all $\mathbf x,\mathbf y\in \varprojlim(C_n\times X,(f_1,\ldots,f_n)_{C_n})$.
\end{definition}
\begin{observation}
Observe that the defined relation $\sim_{(f_1,\ldots,f_n)_{C_n}}$ is an equivalence relation on $\varprojlim(C_n\times X,(f_1,\ldots,f_n)_{C_n})$.
\end{observation}
\begin{theorem}\label{juhuhu}
Let $X$ be a compact metric space, let  $n$ be a positive integer,  let $f_1,f_2,f_3,\ldots,f_n:X\rightarrow X$ be any homeomorphisms, and let  $\sim$ be an equivalence relation on $(C_n\times C_n)\times X$.  Then $\varprojlim(C_n\times X,(f_1,\ldots,f_n)_{C_n})/_{\sim_{(f_1,\ldots,f_n)_{C_n}}}$ is homeomorphic to $(C_n\times C_n)\times X/_{\sim}$.
\end{theorem}
\begin{proof}
Let $\varphi:\varprojlim(C_n\times X,(f_1,\ldots,f_n)_{C_n})/_{\sim_{(f_1,\ldots,f_n)_{C_n}}}\rightarrow (C_n\times C_n)\times X/_{\sim}$ be defined by 
$
\varphi([\mathbf x])=[T(\mathbf x)]
$
for any $\mathbf x\in\varprojlim(C_n\times X,(f_1,\ldots,f_n)_{C_n})$, where $T$ is the homeomorphism from the proof of Theorem \ref{loj}.  Then $\varphi$ is a well-defined function, which is a homeomorphism.
\end{proof}

\begin{theorem}\label{ab}
Let $X$ be a compact metric space, let  $n$ be a positive integer,  and let $f_1,f_2,f_3,\ldots,f_n:X\rightarrow X$ be any homeomorphisms.  If  the $n$-tuple $(f_1, f_2, f_3, \ldots, f_n)$ has property$\mathcal L$, then there is a transitive homeomorphism on the quotient space
$$
\varprojlim(C_n\times X,(f_1,\ldots,f_n)_{C_n})/_{\sim_{(f_1,\ldots,f_n)_{C_n}}}.
$$
\end{theorem}
\begin{proof}
By Theorem \ref{miladojka1}, $\sigma_{(f_1,\ldots,f_n)_{C_n}}$ is a transitive homeomorphism on the inverse limit $\varprojlim(C_n\times X,(f_1,\ldots,f_n)_{C_n})$. Therefore, by Theorem \ref{kvocienti}, the function $\sigma_{(f_1,\ldots,f_n)_{C_n}}^{\star}$ is a transitive homeomorphism on $\varprojlim(C_n\times X,(f_1,\ldots,f_n)_{C_n})/_{\sim_{(f_1,\ldots,f_n)_{C_n}}}$. 
\end{proof}

Next, we construct, using Theorem \ref{juhuhu} and Theorem \ref{ab},  a transitive homeomorphism on the Cantor fan.  First, we prove Lemma \ref{lemmma1} and Lemma \ref{lemmma2}. 

\begin{lemma}\label{lemmma1}
Let $\varepsilon>0$, let $x\in (0,\frac{2}{3}]$ and let $z\in [0,\frac{2}{3}]$. Then there are positive integers $k$ and $n$ such that 
$$
\Big|\Big(\frac{1}{2}\Big)^{\frac{k}{2^n}}\cdot x^{\frac{1}{2^n}}- z\Big| < \varepsilon.
$$ 
\end{lemma}
\begin{proof}
Note that the set $\big\{\frac{k}{2^n} \ | \ k \textup{ and } n \textup{ are positive integers}\big\}$ is dense in $[0,\infty)$. Therefore, $\Big\{\Big(\frac{1}{2}\Big)^{\frac{k}{2^n}} \ | \ k \textup{ and } n \textup{ are positive integers}\Big\}$ is dense in $[0,1]$. For each positive integer $n$, we choose a positive integer $k_n$ as follows.
\begin{enumerate}
\item If $z>0$, then for each positive integer $n$, let $k_n$ be such a positive integer that for each positive integer $k$,
$
\big| (\frac{1}{2})^{\frac{k_n}{2^n}}-z \big|\leq \big| (\frac{1}{2})^{\frac{k}{2^n}}-z \big|.
$
\item If $z=0$, then for each positive integer $n$, let $k_n$ be a positive integer such that  
$
\big(\frac{1}{2}\big)^{\frac{k_n}{2^n}} \leq \frac{1}{2^n}.
$
\end{enumerate}
Then $ \lim_{n\to \infty}\big(\frac{1}{2}\big)^{\frac{k_n}{2^n}}=z$. 
{ Since} $\displaystyle \lim_{n\to \infty}x^{\frac{1}{2^n}}=1$, { it follows that}
$$
\lim_{n\to \infty}\Big(\Big(\frac{1}{2}\Big)^{\frac{k_n}{2^n}}\cdot x^{\frac{1}{2^n}}\Big)=\lim_{n\to \infty}\Big(\frac{1}{2}\Big)^{\frac{k_n}{2^n}}\cdot \lim_{n\to \infty}x^{\frac{1}{2^n}}=z\cdot 1=z 
$$
and {  hence} there is a positive integer $n$ such that 
$
\big|\big(\frac{1}{2}\big)^{\frac{k_n}{2^n}}\cdot x^{\frac{1}{2^n}}- z\big| < \varepsilon.
$
Choose any such $n$ and let $k=k_n$. 
\end{proof}
\begin{definition}\label{definicija}
We use $f_1$, $f_2$ and $f_3$  to denote the homeomorphisms $f_1,f_2,f_3:[0,1]\rightarrow [0,1]$ that are defined by 
$$
f_1(t)=\sqrt x,  ~~~   f_2(t)=\begin{cases}
				\frac{1}{2}t\text{;} & t\leq \frac{2}{3} \\
				2x-1\text{;} & t\geq \frac{2}{3}
			\end{cases} ~~~~  \textup{ and } ~~~~  f_3(t)=f_1^{-1}(t)
$$
for any $t\in [0,1]$. 
\end{definition}
\begin{lemma}\label{lemmma2}
Let $\varepsilon>0$, let $x\in (0,1)$ and let $z\in [0,1]$. Then there are positive integers $k$, $m$ and $n$ such that 
$$
\Big|f_1^{m+n}(f_2^k(f_3^{m}(x)))- z\Big| < \varepsilon.
$$ 
\end{lemma}
\begin{proof}
We consider the following possible cases.
\begin{enumerate}
\item $z=1$. Let $m=k=1$ and let $y=f_1^{m}(f_2^k(f_3^{m}(x)))$. Then $y\in (0,1)$. Also, let $n$ be a positive integer such that $|f_1^n(y)-1|<\varepsilon$. Then
$$
\Big|f_1^{m+n}(f_2^k(f_3^{m}(x)))- z\Big|=\Big|f_1^{n}(y)- z\Big|=\Big|f_1^{n}(y)- 1\Big| < \varepsilon.
$$
\item $z<1$. Let $m$ be a positive integer such that { both} $f_3^{m}(x)<\frac{2}{3}$ and $f_3^{m}(z)<\frac{2}{3}$. 
{Since the map $f_1^m$ is uniformly continuous, there is a $\delta>0$ such that for all $p,q\in [0,1]$,}
$$
|p-q|<\delta ~~~  \Longrightarrow   ~~~   |f_1^m(p)-f_1^m(q)|<\varepsilon.
$$
Choose and fix such a $\delta>0$.  By Lemma  \ref{lemmma1}, there are positive integers $k$ and $n$ such that 
$$
\Big|\Big(\frac{1}{2}\Big)^{\frac{k}{2^n}}\cdot \Big(f_3^{m}(x)\Big)^{\frac{1}{2^n}}- f_3^{m}(z)\Big| < \delta.
$$
Choose and fix such  positive integers $k$ and $n$.  Note that 
$
f_1^{n}(f_2^k(f_3^{m}(x)))=\big(\frac{1}{2}\big)^{\frac{k}{2^n}}\cdot \big(f_3^{m}(x)\big)^{\frac{1}{2^n}}
$
and that
$
z=f_1^m(f_3^{m}(z)).
$
Therefore,
$$
\Big|f_1^{m+n}(f_2^k(f_3^{m}(x)))- z\Big| = \Big|f_1^m\left(\Big(\frac{1}{2}\Big)^{\frac{k}{2^n}}\cdot \Big(f_3^{m}(x)\Big)^{\frac{1}{2^n}}\right)- f_1^m\big(f_3^{m}(z)\big)\Big| < \varepsilon.
$$
\end{enumerate}
\end{proof}
\begin{observation}
Note that it follows from Lemma \ref{lemmma2} that the $3$-tuple $(f_1, f_2,f_3)$ has property$\mathcal L$.
\end{observation}
\begin{definition}
We use $\sim$  to denote the equivalence relation on $(C_3\times C_3)\times [0,1]$, defined by
$$
((\mathbf x_1,\mathbf x_2),t)\sim((\mathbf y_1,\mathbf y_2),s) \Longleftrightarrow  ((\mathbf x_1,\mathbf x_2),t)=((\mathbf y_1,\mathbf y_2),s) \textup{ or } s=t=1
$$
for all $((\mathbf x_1,\mathbf x_2),t),((\mathbf y_1,\mathbf y_2),s)\in (C_3\times C_3)\times [0,1]$.  
\end{definition}
\begin{observation}
Observe that $((C_3\times C_3)\times [0,1])/_{\sim}$ is a Cantor fan.
\end{observation}
Here is our first proof of Theorem \ref{main}.

\begin{proof}{ (Our first proof of Theorem \ref{main})

It follows from Theorem \ref{juhuhu} that 
$
\varprojlim(C_3\times X,(f_1,f_2,f_3)_{C_3})/_{\sim_{(f_1,f_2,f_3)_{C_3}}}
$
is homeomorphic to $((C_3\times C_3)\times [0,1])/_{\sim}$ and it is, therefore, a Cantor fan.  Since $f_1$, $f_2$ and $f_3$ have property $\mathcal L$, it follows from Theorem \ref{ab} that there is a transitive homeomorphism on $\varprojlim(C_3\times X,(f_1,f_2,f_3)_{C_3})/_{\sim_{(f_1,f_2,f_3)_{C_3}}}$.}
\end{proof}

\subsection{The quotient spaces of  $D_n\times X$ }

\begin{theorem}\label{glavD}
Let $X$ be a compact metric space that does not have any isolated points, let $n$ be a positive integer, and let $f_1,f_2,f_3,\ldots,f_n:X\rightarrow X$ be any homeomorphisms.  If the $n$-tuple $(f_1, f_2, f_3, \ldots,f_n)$ has property$\mathcal L$, then  $(f_1,\ldots,f_n)_{D_n}$ is a transitive homeomorphism.
\end{theorem}
\begin{proof}
By Theorem \ref{miladojka1}, the shift homeomorphism $\sigma_{(f_1,\ldots,f_n)_{C_n}}$ on the inverse limit $\varprojlim(C_n\times X,(f_1,\ldots,f_n)_{C_n})$ is transitive. Therefore, $\sigma_{(f_1,\ldots,f_n)_{C_n}}^{-1}$ is transitive. We show that $(f_1,\ldots,f_n)_{D_n}$ is transitive by showing that $(f_1,\ldots,f_n)_{D_n}$ and $\sigma_{(f_1,\ldots,f_n)_{C_n}}^{-1}$ are topological conjugates. Let the homeomorphism $T:\varprojlim(C_n\times X,(f_1,\ldots,f_n)_{C_n})\rightarrow (C_n\times C_n)\times X$ be defined as in the proof of Theorem \ref{loj} by 
$
T(\mathbf x)=((\mathbf a_{\mathbf x},\mathbf b_{\mathbf x}),t_{\mathbf x})
$
for each $\mathbf x\in \varprojlim(C_n\times X,(f_1,\ldots,f_n)_{C_n})$. Also, let $S:(C_n\times C_n)\times X\rightarrow D_n\times X$ be defined by 
$
S((\mathbf a,\mathbf b),t)=(({\ldots, \mathbf b(3),\mathbf b(2),\mathbf b(1)};{\mathbf a(1),\mathbf a(2),\mathbf a(3),\ldots}),t)
$
for each $((\mathbf a,\mathbf b),t)$. Clearly, $S$ is a homeomorphism. 
One can easily check that 
$$
(f_1,\ldots,f_n)_{D_n}(\mathbf x,t)=(S\circ T\circ \sigma_{(f_1,\ldots,f_n)_{C_n}}^{-1}\circ T^{-1}\circ S^{-1})(\mathbf x,t)
$$
for each $(\mathbf x,t)\in D_n\times X$.
\end{proof}

\begin{theorem}\label{taprvi}
Let $X$ be a compact metric space, let  $n$ be a positive integer,  and let $f_1,f_2,f_3,\ldots,f_n:X\rightarrow X$ be any homeomorphisms.  If  $f_1$, $f_2$, $f_3$, $\ldots$,  $f_n$ have property $\mathcal L$, then there is a transitive homeomorphism on $D_n\times X/_{\sim}$ for any equivalence relation $\sim$ on $D_n\times X$.
\end{theorem}
\begin{proof}
Let $\sim$ be an equivalence relation on $D_n\times X$. By Theorem \ref{glavD}, there is a transitive homeomorphism on $D_n\times X$. Let $h:D_n\times X\rightarrow D_n\times X$ be such a homeomorphism.  By Theorem \ref{kvocienti}, $h^{\star}$ is a transitive homeomorphism on $D_n\times X/_{\sim}$.
\end{proof}
\begin{definition}
We use $\sim$ to denote the equivalence relation on $D_3\times [0,1]$, defined by
$$
(\mathbf x,t)\sim(\mathbf y,s) \Longleftrightarrow  (\mathbf x,t)=(\mathbf y,s) \textup{ or } s=t=1
$$
for all $(\mathbf x,t),(\mathbf y,s)\in D_3\times [0,1]$.  
\end{definition}
\begin{observation}
Observe that $(D_3\times [0,1])/_{\sim}$ is a Cantor fan.
\end{observation}
We conclude this section by giving our second proof  of Theorem \ref{main}.

\begin{proof}{ (Our second proof of Theorem \ref{main})

Let $f_1$, $f_2$ and $f_3$ be the homeomorphisms from Definition \ref{definicija}.  Recall that $f_1$, $f_2$ and $f_3$ have property$\mathcal L$.  Therefore,  by Theorem \ref{taprvi},  there is a transitive homeomorphism on $(D_3\times [0,1])/_{\sim}$. }
\end{proof}
\section{Transitive homeomorphisms on Mahavier products of closed relations}\label{s3}
Here, we study {relationships} between  $X_F$ and $X_F^+$. In Subsection \ref{unka},  we show that the shift homeomorphism on $X_F$ is transitive if and only if the shift map on $X_F^+$ is transitive.  In Subsection \ref{dunka},  we study closed relations on $X$ that are unions of graphs of continuous functions.  Then, we apply our results to obtain transitive homeomorphisms on the Cantor fan.
\subsection{Two-sided Mahavier products}\label{unka}
Theorems \ref{povezava} and \ref{tazadnji} are the main results of Subsection \ref{unka},  where a relation between $X_F$ and $X_F^+$ is established. 
\begin{theorem}\label{povezava}
Let $X$ be a compact metric space and let $F$ be a closed relation on $X$. Then $\varprojlim(X_F^{+},\sigma_F^+)$ is homeomorphic to the two-sided Mahavier product $X_F$. Also, the inverse of the shift map $\sigma_F$ on $X_F$ is topologically conjugate to the  shift map on $\varprojlim(X_F^{+},\sigma_F^+)$. 
\end{theorem}
\begin{proof}
Note that for any $\mathbf x\in \varprojlim(X_F^{+},\sigma_F^+)$, there are unique points $\mathbf a\in X_F^+$ and $\mathbf b\in X_{F^{-1}}^+$ such that
$
\mathbf x=(\mathbf a,(\mathbf b(1))\oplus \mathbf a, (\mathbf b(2),\mathbf b(1))\oplus \mathbf a, (\mathbf b(3),\mathbf b(2),\mathbf b(1))\oplus \mathbf a, \ldots).
$ 
Note that in this case, $(\mathbf b(1),\mathbf a(1))\in F$. 
So, for each $\mathbf x\in \varprojlim(X_F^{+},\sigma_F^+)$, let $\mathbf a_{\mathbf x}\in X_F^+$ and $\mathbf b_{\mathbf x}\in X_{F^{-1}}^+$ be such that
$
\mathbf x=(\mathbf a_{\mathbf x},(\mathbf b_{\mathbf x}(1))\oplus \mathbf a_{\mathbf x}, (\mathbf b_{\mathbf x}(2),\mathbf b_{\mathbf x}(1))\oplus \mathbf a_{\mathbf x}, (\mathbf b_{\mathbf x}(3),\mathbf b_{\mathbf x}(2),\mathbf b_{\mathbf x}(1))\oplus \mathbf a_{\mathbf x}, \ldots).
$ 
Let $\varphi:\varprojlim(X_F^{+},\sigma_F^+)\rightarrow X_{F}$ be defined by 
$$
\varphi(\mathbf x)=({\ldots, \mathbf b_{\mathbf x}(3),\mathbf b_{\mathbf x}(2),\mathbf b_{\mathbf x}(1)}; {\mathbf a_{\mathbf x}(1),\mathbf a_{\mathbf x}(2),\mathbf a_{\mathbf x}(3),\ldots})
$$
for any $\mathbf x\in \varprojlim(X_F^{+},\sigma_F^+)$.  Then $\varphi$ is a homeomorphism.  Next, let $\sigma$ be the shift map on $\varprojlim(X_F^{+},\sigma_F^+)$.  Then one can easily check that $\sigma_F^{-1}=\varphi\circ \sigma \circ \varphi^{-1}$.
\end{proof}
Next, we give an example of {a} dynamical system $(X,f)$ such that $(\varprojlim(X,f),\sigma)$ is transitive while $(X,f)$ is not ($\sigma :\varprojlim(X,f)\rightarrow \varprojlim(X,f)$ is the shift map  on $\varprojlim(X,f)$).
\begin{example}
Let $f:[0,1]\rightarrow [0,1]$ be defined by $f(t)=0$ for any $t\in [0,1]$. Then $\varprojlim([0,1],f)=\{(0,0,0,\ldots)\}$ and it follows that $(\varprojlim([0,1],f),\sigma)$ is transitive.  Since $f$ is not surjective, $( {[0,1]},f)$ is not a transitive dynamical system.
\end{example}
However, the following holds.
\begin{proposition}\label{labelo}
Let $X$ be a compact metric space, let $f:X\rightarrow X$ be a continuous surjection and let $\sigma :\varprojlim(X,f)\rightarrow \varprojlim(X,f)$ be the shift map on $\varprojlim(X,f)$. The following statements are equivalent.
\begin{enumerate}
\item\label{unc}  $(X,f)$ is transitive.
\item\label{uncc}  $(\varprojlim(X,f),\sigma)$ is transitive.
\end{enumerate}
\end{proposition}
\begin{proof}
The implication from \ref{unc} to \ref{uncc} follows from Theorem \ref{shifttransitive}. To prove the implication from \ref{uncc} to \ref{unc}, suppose that $(\varprojlim(X,f),\sigma)$ is transitive. By Observation \ref{isolatedpoints}, also $(\varprojlim(X,f),\sigma^{-1})$ is transitive since $\sigma$ is a homeomorphism.  Let $U$ and $V$ be any non-empty open sets in $X$, and let $W_U=U\times \prod_{k=2}^{\infty}X$ and $W_V=V\times \prod_{k=2}^{\infty}X$. Also, let $\mathbf x\in W_U\cap \varprojlim(X,f)$. Since $(\varprojlim(X,f),\sigma^{-1})$ is transitive, there is a non-negative integer $n$ such that $(\sigma^{-1})^n(\mathbf x)\in W_V\cap \varprojlim(X,f)$. It follows that for such an $n$, $f^n(\mathbf x(1))\in V$, therefore,  $f^n(\mathbf x(1))\in f^n(U)\cap V$ and $f^n(U)\cap V\neq \emptyset$ follows.
\end{proof}
\begin{theorem}\label{tazadnji}
Let $X$ be a compact metric space and let $F$ be a closed relation on $X$. The following statements are equivalent. 
\begin{enumerate}
\item The map $\sigma_F^+$ is transitive.
\item The homeomorphism $\sigma_F$ is transitive. 
\end{enumerate}
\end{theorem}
\begin{proof}
Let $\sigma$ be the shift map on $\varprojlim(X_F^{+},\sigma_F^+)$. First, suppose that $\sigma_F^+$ is transitive.  By  Theorem \ref{shifttransitive}, $\sigma$ is also transitive.  By Theorem \ref{povezava},  {$\sigma$} is topologically conjugate to $\sigma_F^{-1}$, therefore, $\sigma_F^{-1}$ is transitive. It follows that $\sigma_F$ is transitive. 

Now, suppose that $\sigma_F$ is transitive. Since $\sigma_F^{-1}$ is transitive, it follows from Theorem \ref{povezava} that $\sigma$ is transitive. Since $\sigma_F^+$ is surjective, it follows from Proposition \ref{labelo} that $\sigma_F^+$ is transitive.
\end{proof}
\subsection{Unions of graphs of continuous functions}\label{dunka}
Theorems \ref{main21} and \ref{main22} are the main results of this sections.  After proving them, we apply them to obtain transitive homeomorphisms on the Cantor fan. First, we prove the following lemma. 
\begin{lemma}\label{L1}
Let $X$ be a compact metric space and let $\{f_{\alpha} \ | \ \alpha \in \Lambda\}$ be a collection of continuous functions from $X$ to $X$ such that $F=\bigcup_{\alpha\in \Lambda}\Gamma(f_{\alpha})$ is closed in $X\times X$.  Then for each open set $U$ in $\prod_{k=1}^{\infty}X$ and for each $\mathbf x\in U\cap X_F^+$,  
$$
\mathbf x(1)\in \Int (\pi_1(U\cap X_F^+)).
$$ 
\end{lemma}
\begin{proof}
Let $U$ be an open set in $\prod_{k=1}^{\infty}X$ such that $U\cap X_F^+\neq \emptyset$ and let $\mathbf x\in U\cap X_F^+$.  Also, let $n$ be a positive integer and let $U_1$, $U_2$, $U_3$, $\ldots$,  $U_n$ be open sets in $X$ such that 
$$
\mathbf x\in \Big(U_1\times U_2\times U_3\times \ldots \times U_{n-1}\times U_n\times \prod_{k=n+1}^{\infty}X\Big)\cap X_F^+\subseteq U\cap X_F^+.
$$
Next, let $\alpha_1,\alpha_2,\alpha_3,\ldots,\alpha_{n-1}\in \Lambda$ be such that $\mathbf x(k+1)=f_{\alpha_k}\big(\mathbf x(k)\big)$ for each $k\in \{1,2,3,\ldots, n-1\}$ and let $V_n=U_n$.   Since $f_{\alpha_{n-1}}$ is continuous at the point $\mathbf x(n-1)$,  there is an open set $V_{n-1}$ in $X$ such that $\mathbf x(n-1)\in V_{n-1}\subseteq U_{n-1}$ and $f_{\alpha_{n-1}}(V_{n-1})\subseteq  {V_n}$.  Since $f_{\alpha_{n-2}}$ is continuous at the point $\mathbf x(n-2)$,  there is an open set $V_{n-2}$ in $X$ such that $\mathbf x(n-2)\in V_{n-2}\subseteq U_{n-2}$ and $f_{\alpha_{n-2}}(V_{n-2})\subseteq V_{n-1}$.  Continuing inductively,  we construct open sets $V_1$, $V_2$, $V_3$, $\ldots$,   $V_n$ in $X$ such that for each $k\in \{1,2,3,\ldots, n-1\}$,  $\mathbf x(k)\in V_k\subseteq U_k $ and $f_{\alpha_{k}}(V_{k})\subseteq V_{k+1}$.

 {Finally, w}e show that $\mathbf x(1)\in \Int (\pi_1(U\cap X_F^+))$.  Since $\mathbf x(1)\in V_1$ and $V_1$ is open in $X$,  it suffices to see that $V_1\subseteq \pi_1(U\cap X_F^+)$.  Let $x\in V_1$.  To see that $x\in \pi_1(U\cap X_F^+)$, let $y_1=x$,  for each $k\in \{2,3,\ldots,n\}$, let $y_k=f_{\alpha_{k-1}}(x_{k-1})$, for each $k>n$, let $y_k=f_{\alpha_{1}}(x_{k-1})$, and let $\mathbf y=(y_1,y_2,y_3,\ldots)$.  Then $\mathbf y\in U\cap X_F^+$ and, therefore, $\mathbf y(1)\in \pi_1(U\cap X_F^+)$. It follows that $x\in  \pi_1(U\cap X_F^+)$.
\end{proof}
The following observation follows directly from the proof of Lemma \ref{L1}.
\begin{observation}\label{O1}
Let $X$ be a compact metric space,  let $\{f_{\alpha} \ | \ \alpha \in \Lambda\}$ be a collection of continuous functions from $X$ to $X$ such that $F=\bigcup_{\alpha\in \Lambda}\Gamma(f_{\alpha})$ is closed in $X\times X$,   let $U$ be an open set in $\prod_{k=1}^{\infty}X$,  let $\mathbf x\in U\cap X_F^+$,  and for each positive integer $k$,   {let} $\alpha_k\in \Lambda$ be such that $\mathbf x(k+1)=f_{\alpha_k}\big(\mathbf x(k)\big)$.  Then there are a positive integer $m$ and open sets $U_1$, $U_2$, $U_3$, $\ldots$,  $U_m$ in $X$ such that 
\begin{enumerate}
\item $\mathbf x\in \Big(U_1\times U_2\times U_3\times \ldots \times U_{m-1}\times U_m\times \prod_{k=m+1}^{\infty}X\Big)\cap X_F^+\subseteq U\cap X_F^+$ and 
\item for each $k\in \{1,2,3,\ldots, m-1\}$,  
$
f_{\alpha_{k}}(U_{k})\subseteq U_{k+1}.
$
\end{enumerate}
\end{observation}
\begin{definition}
Let $X$ be a compact metric space,  let $F$ be a closed relation on $X$ and let $x\in X$. Then we define 
$$
\mathcal U^{\oplus}_F(x)=\{y\in X \ | \ \textup{there are } n\in \mathbb N \textup{ and } \mathbf x\in X_F^{n} \textup{ such that } \mathbf x(1)=x \textup{ and } \mathbf x(n)=y \}
$$
and we call it the forward impression of $x$ by $F$.
\end{definition}

\begin{theorem}\label{main21}
{ Let $X$ be a compact metric space, let $F$ be a closed relation on $X$,  let  $\{f_{\alpha} \ | \ \alpha \in A\}$ be a non-empty collection of continuous functions from $X$ to $X$ such that $F^{-1}=\bigcup_{\alpha\in A}\Gamma(f_{\alpha})$,  and let  $\{g_{\beta} \ | \ \beta \in B\}$ be a non-empty collection of continuous functions from $X$ to $X$ such that $F=\bigcup_{\beta\in B}\Gamma(g_{\beta})$. } If there is a dense set $D$ in $X$ such that for each $s\in D$,  $\Cl(\mathcal U^{\oplus}_F(s))=X$, then $(X_F^+,\sigma_F^+)$ is transitive. 
\end{theorem}
\begin{proof}
Let $U$ and $V$ be open sets in $\prod_{k=1}^{\infty}X$ such that $U\cap X_F^+\neq \emptyset$ and $V\cap X_F^+\neq \emptyset$,  let  $\mathbf x\in U\cap X_F^+$,  and let for each positive integer $k$, {
\begin{enumerate}
\item $\alpha_k\in A$ be such that $\mathbf x(k)=f_{\alpha_k}\big(\mathbf x(k+1)\big)$, and
\item $\beta_k\in B$ be such that $\mathbf x(k+1)=g_{\beta_k}\big(\mathbf x(k)\big)$. 
\end{enumerate}
Also, let $m$ be a positive integer and let $U_1$, $U_2$, $U_3$, $\ldots$,  $U_m$ be non-empty open sets in $X$ such that 
\begin{enumerate}
\item $\mathbf x\in \Big(U_1\times U_2\times U_3\times \ldots \times U_{m-1}\times U_m\times \prod_{k=m+1}^{\infty}X\Big)\cap X_F^+\subseteq U\cap X_F^+$ and 
\item for each $k\in \{1,2,3,\ldots, m-1\}$,  
$
g_{\beta_{k}}(U_{k})\subseteq U_{k+1}.
$
\end{enumerate}
 Such an integer and such sets do exist by Observation \ref{O1}. 
By Lemma \ref{L1}, $\mathbf x(1)\in\Int (\pi_1(U\cap X_F^+))$.  Let $V_1 = U_1\cap \Int (\pi_1(U\cap X_F^+))$. Then $V_1$ is open in $X$, $V_1\subseteq U_1$ and $\mathbf x(1)\in V_1$.  Next, let $V_2$ be an open set in $X$ such that $V_2\subseteq U_2$,  $\mathbf x(2)\in V_2$ and $f_{\alpha_1}(V_2)\subseteq V_1$.  Since $f_{\alpha_1}:X\rightarrow X$ is continuous, it follows that such a set $V_2$ does exist.  Continuing inductively, we construct open sets $V_1$, $V_2$, $V_3$, $\ldots$,  $V_m$ in $X$ such that for each $k\in \{1,2,3,\ldots, m-1,m\}$,  $V_k\subseteq U_k$, $\mathbf x(k)\in V_k$ and for each $k\in \{1,2,3,\ldots, m-1\}$,  $f_{\alpha_k}(V_{k+1})\subseteq V_{k}$. 

Next, let $d_m\in V_m\cap D$, $d_{m-1}=f_{\alpha_{m-1}}(d_m)$, $d_{m-2}=f_{\alpha_{m-2}}(d_{-1})$, $\ldots$, $d_{2}=f_{\alpha_{2}}(d_3)$ and $d_{1}=f_{\alpha_{1}}(d_2)$. Also, let $\mathbf d\in X_F^+$ be such a point that $\mathbf d( {1})=g_{\beta_m}(d_m)$ and let 
$$
\mathbf x'=(d_1,d_2,d_3,\ldots,d_m)\oplus \mathbf d.
$$}
   Note that $\mathbf x'\in \Big(U_1\times U_2\times U_3\times \ldots \times U_{m-1}\times U_m\times \prod_{k=m+1}^{\infty}X\Big)\cap X_F^+$ is such a point that $\mathbf x'(m)\in U_m\cap D$. Also, note that $\mathbf x'\in U\cap X_F^+$. Since for each $s\in D$,  $\Cl(\mathcal U^{\oplus}_F(s))=X$,  it follows that $\Cl(\mathcal U^{\oplus}_F(\mathbf x'(m)))=X$.  
   
   It follows from Lemma \ref{L1} that $\Int (\pi_1(V\cap X_F^+))\neq \emptyset$.  Therefore, it follows from $\Cl(\mathcal U^{\oplus}_F(\mathbf x'(m)))=X$ that there are a positive integer $\ell$ and a point $\mathbf y\in X_F^{\ell}$ such that $\mathbf y(1)=\mathbf x'(m)$ and $\mathbf y(\ell)\in \Int (\pi_1(V\cap X_F^+))$. Finally, let $\mathbf z\in V\cap X_F^+$ be such that $\mathbf z(1)=\mathbf y(\ell)$ and let
$$
\mathbf w=(\mathbf x'(1),\mathbf x'(2),\mathbf x'(3),\ldots,\mathbf x'(m),\mathbf y(2),\mathbf y(3),\mathbf y(4),\ldots,\mathbf y(\ell),\mathbf z(2),\mathbf z(3),\mathbf z(4),\ldots).
$$
Then $\mathbf w\in U\cap X_F^+$ and $(\sigma_F^+)^{m+\ell}(\mathbf w)=\mathbf z\in V\cap X_F^+$.  Therefore, $(\sigma_F^+)^{m+\ell}(U\cap X_F^+)\cap (V\cap X_F^+)\neq \emptyset$ and it follows that $(X_F^+,\sigma_F^+)$ is transitive. 
\end{proof}
The following well-known result, which is presented in Lemma \ref{lemissimamissima}, {is used later in the paper.  }
\begin{lemma}\label{lemissimamissima}
For each $x\in (0,1)$, the set $\{x^{\frac{2^m}{3^n}} \ | \  m \textup{ and } n \textup{ are positive integers}\}$ is dense in $[0,1]$.
\end{lemma}
\begin{proof}
An irrational rotation on a circle produces a dense orbit of any point on a circle. This translates to the fact that for any irrational number $\gamma$, the set $\{m-n\cdot \gamma \ | \ m \textup{ and } n \textup{ are positive integers}\}$ is dense in $\mathbb R$. Therefore, the set $\{m-n\cdot \frac{\ln 3}{\ln 2} \ | \ m \textup{ and } n \textup{ are positive integers}\}$ is dense in $\mathbb R$, since $\frac{\ln 3}{\ln 2}$ is irrational. Thus, the set $\{m\cdot \ln 2-n\cdot \ln 3 \ | \ m \textup{ and } n \textup{ are positive integers}\}$ is dense in $\mathbb R$. Therefore, the set $\{e^{m\cdot \ln 2-n\cdot \ln 3 } \ | \ m \textup{ and } n \textup{ are positive integers}\}$ is dense in $[0,\infty)$ and it follows that the set 
$\{\frac{2^m}{3^n} \ | \ m \textup{ and } n \textup{ are positive integers}\}$ is dense in $[0,\infty)$. 

Let $x\in (0,1)$ and let $a,b\in (0,1)$ be such that $a<b$. Then $\ln a<\ln b<0$ and $\ln x<0$. So $0<\frac{\ln b}{\ln x}<\frac{\ln a}{\ln x}$. 
Let $m$ and $n$ be positive integers such that 
$$
\frac{\ln b}{\ln x}<\frac{2^m}{3^n}<\frac{\ln a}{\ln x}. 
$$
Such integers $m$ and $n$ do exist since $\{\frac{2^m}{3^n} \ | \ m \textup{ and } n \textup{ are positive integers}\}$ is dense in $[0,\infty)$.  Then $\ln a<\frac{2^m}{3^n}\cdot \ln x<\ln b$ and, therefore, $\ln a<\ln x^{\frac{2^m}{3^n}}<\ln b$. It follows that $a<x^{\frac{2^m}{3^n}}<b$.  
\end{proof}
In the following example, we show that in general, there are closed relations $F$ on compact metric spaces $X$ such that $(X_F^+,\sigma_F^+)$ is not transitive even though there is a dense set $D$ in $X$ such that for each $s\in D$,  $\Cl(\mathcal U^{\oplus}_F(s))=X$.
\begin{example}
Let 
$$
F=\{(x,y)\in [0,1]\times [0,1] \ | \ y=x^2 \textup{ or } y=x^{\frac{1}{3}}\}\cup\{(0,1),(1,0)\}.
$$
Then $F$ is a closed relation on $[0,1]$ such that 
\begin{enumerate}
\item for each $x\in (0,1)$,  $\mathcal U^{\oplus}_F(x)$ is dense in $[0,1]$ {(by Lemma \ref{lemissimamissima})},
\item for each $x\in (0,1)$,  $\mathcal U^{\oplus}_F(x)\cap\{0,1\}=\emptyset$.
\end{enumerate}
Let $U=\big((0,1)\times \prod_{k=2}^{\infty}[0,1]\big)\cap X_F^+$ and $V=\big([0,\frac{1}{100})\times (\frac{99}{100},1]\times \prod_{k=3}^{\infty}[0,1]\big)\cap X_F^+$. Then $U\neq \emptyset$, $V\neq \emptyset$ and $\pi_1(V)=\{0\}$. Therefore,  for each non-negative integer $n$,  $(\sigma_F^+)^n(U)\cap V=\emptyset$ and it follows that the dynamical system  $([0,1]_F^+, \sigma_F^+)$ is not transitive.
\end{example}

In the following example, we demonstrate how Theorem \ref{main21} may be used to obtain a well-known result. 
\begin{example}
Let $f:[0,1]\rightarrow [0,1]$ be the tent-map defined by $f(x)=2x$, if $x\leq \frac{1}{2}$, and $f(x)=2-2x$, if $x\geq \frac{1}{2}$.  Also, let $F=\Gamma(f)^{-1}$.  {It is a well-known fact that the inverse limit of closed intervals $[0,1]$ with a tent map as a bonding function is a Knaster continuum. Since   $[0,1]_F^+$ is the inverse limit of closed intervals with a tent map as a bonding function, it is the Knaster continuum.  It is also a well-known fact that the shift homeomorphism on $[0,1]_F^+$ is transitive. }  Note that this  follows easily also from Theorem \ref{main21} since 
\begin{enumerate}
\item $F^{-1}$ is the graph of $f$,
\item $F$ is the union of graphs of two continuous functions,  and 
\item there is a dense subset of $[0,1]$ with a dense forward impression. 
\end{enumerate} 
\end{example}

\begin{definition}
Let $X$ be a compact metric space (with metric $d$) and let $F$ be a closed relation on $X$.  We use  $D_F$ to denote the metric on $X_F$ and $D_F^{+}$ to denote the metric on $X_F^{+}$:
$$
D_F(\mathbf x,\mathbf y)=\sum_{k=-\infty}^{\infty}\frac{d(\mathbf x(k),\mathbf y(k))}{2^{|k|}}  ~~~  \textup{ and }  ~~~    D_F^+(\mathbf x,\mathbf y)=\sum_{k=1}^{\infty}\frac{d(\mathbf x(k),\mathbf y(k))}{2^k}
$$
for all $\mathbf x,\mathbf y\in X_F$ and $\mathbf x,\mathbf y\in X_F^+$.
\end{definition}
\begin{theorem}\label{main22}
Let $X$ be a compact metric space, let $n$ be a positive integer,  and for each $k\in \{1,2,3,\ldots,n\}$, let $f_k:X\rightarrow X$ be a homeomorphism. Let $F=\bigcup_{k=1}^n\Gamma(f_k)$  and assume that for each $x\in X$ and for all $i,j\in \{1,2,3,\ldots,n\}$,
$$
f_i(x)=f_j(x) \textup{ and } i\neq j ~~~ \Longrightarrow ~~~   F(x)=\{x\}.
$$
Then $X_F$ is homeomorphic to $X_F^+$.
\end{theorem}
\begin{proof}
Without any loss of generality, we assume that the diameter of $X$ is less or equal to 1. Let 
$
A=\{a\in X \ | \ \textup{ for each } k\in \{1,2,3,\ldots ,n\}, f_k(a)=a\}
$
and let 
$
A^{\star}=\{(\ldots,a,a;a,a,a,\ldots) \ | \ a\in A\}.
$
Note that $A$ is a closed subset of $X$ and that $A^{\star}$ is a closed subset of $X_F$. Also, note that for each $\mathbf x\in X_F\setminus A^{\star}$, there is a unique function $m_{\mathbf x}:\mathbb Z\rightarrow \{1,2,3,\ldots ,n\}$ such that for each $k\in \mathbb Z$, $\mathbf x (k+1)=f_{m_{\mathbf x}(k)}(\mathbf x(k))$.  For each $\mathbf x\in X_F$, we define $x_1=f_{m_{\mathbf x}(0)}(\mathbf x(0))$, $x_2=f_{m_{\mathbf x}(-1)}(x_1)$, $x_3=f_{m_{\mathbf x}(1)}(x_2)$,  $x_4=f_{m_{\mathbf x}(-2)}(x_3)$, and for each positive integer $k$,  we define
$$
x_{k+1}=\begin{cases}
				f_{m_{\mathbf x}(\frac{k}{2})}(x_k) \text{;} & k \textup{ is even}\\
				f_{m_{\mathbf x}(-\frac{k+1}{2})}(x_k) \text{;} & k \textup{ is odd}.
			\end{cases}
$$
Then we define the function $T:X_f\rightarrow X_F^+$ by
$$
T(\mathbf x)=\begin{cases}
				(x_1,x_2,x_3,\ldots)\text{;} & \mathbf x\in X_F\setminus A^{\star} \\
				(\mathbf x(1),\mathbf x(2),\mathbf x(3),\ldots) \text{;} & \mathbf x\in A^{\star}.
			\end{cases}
$$
Let $\mathbf x\in X_F$ be any point. We show that $T$ is continuous at the point $\mathbf x$. Let $\varepsilon >0$. We consider the following two possible cases.
\begin{enumerate}
\item $\mathbf x\in X_F\setminus A^{\star}$. Let $m$ be a positive integer such that $\sum_{k=m}^{\infty}\frac{1}{2^k}<\frac{\varepsilon}{3}$.  Note that for each $i\in \{-m,-m+1,-m+2,\ldots,m-2,m-1,m\}$ and for each $k\in \{1,2,3,\ldots,n\}\setminus \{m_{\mathbf x}(i)\}$, $f_k(\mathbf x(i))\neq \mathbf x(i+1)$. Let 
$$
U=\Big(\prod_{k=-\infty}^{-m-1}X\Big)\times \Big(\prod_{k=-m}^{m}U_k\Big) \times \Big(\prod_{k=m+1}^{\infty}X\Big),
$$
where for each $k\in \{-m,-m+1,-m+2,\ldots,m-2,m-1,m\}$, $U_k$ is open in $X$, such that $U\cap A^{\star}=\emptyset$ and such that for each $\mathbf y\in U\cap X_F$,   
$\mathbf y(i+1)=f_{m_{\mathbf x}(i)}(\mathbf y(i))$ for each $i\in \{-m,-m+1,-m+2,\ldots,m-2,m-1\}$.  Let  $T(\mathbf x)=(x_1,x_2,x_3,\ldots)$. Then $x_1=f_{m_{\mathbf x}(0)}(\mathbf x(0))$ and  for each positive integer $k$, 
$$
x_{k+1}=\begin{cases}
				f_{m_{\mathbf x}(\frac{k}{2})}(x_k) \text{;} & k \textup{ is even}\\
				f_{m_{\mathbf x}(-\frac{k+1}{2})}(x_k) \text{;} & k \textup{ is odd}.
			\end{cases}
$$
Also,  for each $\mathbf y\in U\cap X_F$, if $T(\mathbf y)=(y_1,y_2,y_3,\ldots)$, then $y_1=f_{m_{\mathbf x}(0)}(\mathbf y(0))$ and  for each positive integer $k\in \{1,2,3,\ldots ,2m+2\}$, 
$$
y_{k+1}=\begin{cases}
				f_{m_{\mathbf x}(\frac{k}{2})}(y_k) \text{;} & k \textup{ is even}\\
				f_{m_{\mathbf x}(-\frac{k+1}{2})}(y_k) \text{;} & k \textup{ is odd}.
			\end{cases}
$$
Let $\delta_1>0$ be such that for each $\mathbf y\in U\cap X_F$,
$$
d(\mathbf x(0),\mathbf y(0))<\delta_1 ~~~  \Longrightarrow ~~~  \sum_{k=1}^{2m+2}\frac{d(x_k,y_k)}{2^k}<\frac{\varepsilon}{3},
$$
where $T(\mathbf y)=(y_1,y_2,y_3,\ldots)$.  Such a $\delta_1$ does exist since all the functions $f_1$, $f_2$, $f_3$, $\ldots$,  $f_m$ are continuous.  Also, let $\delta>0$ be such that for all $\mathbf y\in X_F$, 
$$
D_F(\mathbf x,\mathbf y)<\delta ~~~   \Longrightarrow ~~~  \mathbf y\in U\cap X_F \textup{ and } d(\mathbf x(0),\mathbf y(0))<\delta_1.
$$
Therefore, for each $\mathbf y\in X_F$ such that $D_F(\mathbf x,\mathbf y)<\delta$, 
$$
D_F^+(T(\mathbf x),T(\mathbf y))\leq \sum_{k=1}^{2m+2}\frac{d(x_k,y_k)}{2^k}+\sum_{k=2m+3}^{\infty}\frac{1}{2^k}<\frac{\varepsilon}{3}+\frac{\varepsilon}{3}<\varepsilon,
$$
where $T(\mathbf y)=(y_1,y_2, {y}_3,\ldots)$.  
\item $\mathbf x\in A^{\star}$. Let $m$ be a positive integer such that $\sum_{k=m}^{\infty}\frac{1}{2^k}<\frac{\varepsilon}{3}$ and let $a\in A$ be such that $\mathbf x=(\ldots,a,a;a,a,a,\ldots)$.  Since all the functions $f_1$, $f_2$, $f_3$, $\ldots$,  $f_m$ are continuous and since for each $k\in\{1,2,3,\ldots, m\}$,  $f_k(a)=a$, there is a $\delta_1>0$ such that for each $\mathbf y\in X_F^+$,
$$
d(\mathbf y(0),a)<\delta_1 ~~~ \Longrightarrow ~~~   \sum_{k=1}^{m-1}\frac{d(\mathbf y(k),a)}{2^k}<\frac{\varepsilon}{3}.
$$
Choose and fix such a $\delta_1$. Also, let $\delta>0$ be such that for each $\mathbf y\in X_F$,
$$
D_F(\mathbf x,\mathbf y)<\delta ~~~  \Longrightarrow  ~~~  d(\mathbf y(0),a)<\delta_1.   
$$
Therefore, for each $\mathbf y\in X_F$ such that $D_F(\mathbf x,\mathbf y)<\delta$, 
$$
D_F^+(T(\mathbf x),T(\mathbf y))=D_F^+((a,a,a,\ldots),T(\mathbf y))\leq \sum_{k=1}^{m-1}\frac{d(a,y_k)}{2^k}+\sum_{k=m}^{\infty}\frac{1}{2^k}<\frac{\varepsilon}{3}+\frac{\varepsilon}{3}<\varepsilon,
$$
where $T(\mathbf y)=(y_1,y_2, {y}_3,\ldots)$. 
\end{enumerate}
We have just proved that $T$ is continuous.  Note that since all the functions $f_1$, $f_2$, $f_3$, $\ldots$,  $f_m$ are surjective, also $T$ is surjective.  Finally, we prove that $T$ is injective.  Let $\mathbf x,\mathbf y\in X_F$ be such that $T(\mathbf x)=T(\mathbf y)$ and let 
$$
T(\mathbf x)=(x,g_1(x),g_2(g_1(x)),\ldots) \textup{ and } T(\mathbf y)=(y,h_1(x),h_2(h_1(x)),\ldots),
$$
 where for each positive integer $k$, $g_k,h_k\in \{f_1,f_2,f_3,\ldots,f_n\}$.  It follows that $x=y$.  First, suppose that $g_1(x)=x$. Then $\mathbf x\in A^{\star}$.  Since $h_1(y)=g_1(x)$, it follows that $h_1(y)=y$ and, therefore, $\mathbf y\in A^{\star}$. Therefore, $\mathbf x=\mathbf y$.  Next, suppose that $g_1(x)\neq x$.  Then $h_1(y)\neq y$ and it follows that $\mathbf x,\mathbf y\in X_F\setminus  A^{\star}$.  Since for each $x\in X$ and for all $i,j\in \{1,2,3,\ldots,n\}$,
$$
f_i(x)=f_j(x) \textup{ and } i\neq j ~~~ \Longrightarrow ~~~   F(x)=\{x\},
$$
it follows that for each positive integer $k$, $g_k=h_k$. Therefore, $\mathbf x=\mathbf y$.  We have just proved that $T$ is a continuous bijection from the compactum $X_F$ to the metric space $X_F^+$. Therefore, it is a homeomorphism. 
\end{proof}

We conclude this section by stating various examples. 
\begin{example}\label{exx1}
Let $X=[0,1]\cup [2,3]$ and let 
$
F=\{(x,x^2) \ | \ x\in [0,1]\}\cup \{(x,x^{\frac{1}{3}}) \ | \  x\in [0,1]\}.
$
%
Note that $X_F^+$ is the topological suspension of the Cantor set (see \cite[page 42]{nadler} for the definition of the topological suspension).  By Theorem \ref{main22}, $X_F$ is homeomorphic to $X_F^+$ and by Theorem \ref{main21}, $\sigma_F^+$ is transitive, since for $D=(0,1)$,  the set $\mathcal U^{\oplus}_F(s)$ is dense in $X$ (by Lemma \ref{lemissimamissima}) for each $s\in D$.  It follows from Theorem \ref{tazadnji} that $\sigma_F$ is a transitive homeomorphism on $X_F$.  
\end{example}
In the following examples, transitive homeomorphisms on the Cantor fan are obtained. These are our third and fourth variant of the proof of Theorem \ref{main}.

\begin{example}\label{exx2}{ (Our third proof of Theorem \ref{main})

Let $X=[0,1]{\cup [2,3]}$, let 
$$
f_1(x)=\begin{cases}
				x^2\text{;} & x\in [0,1]\\
				(x-2)^{\frac{1}{3}}+2\text{;} & x\in [2,3]
			\end{cases}     ~~~  \textup{ and }  ~~~  f_2(x)=\begin{cases}
				x+2\text{;} & x\in [0,1]\\
				x-2\text{;} & x\in [2,3]
			\end{cases} 
$$
for each $x\in X$,  and let 
$
F=\Gamma(f_1)\cup \Gamma(f_2);
$
 see Figure \ref{fig2}.
\begin{figure}[h!]
	\centering
		\includegraphics[width=15em]{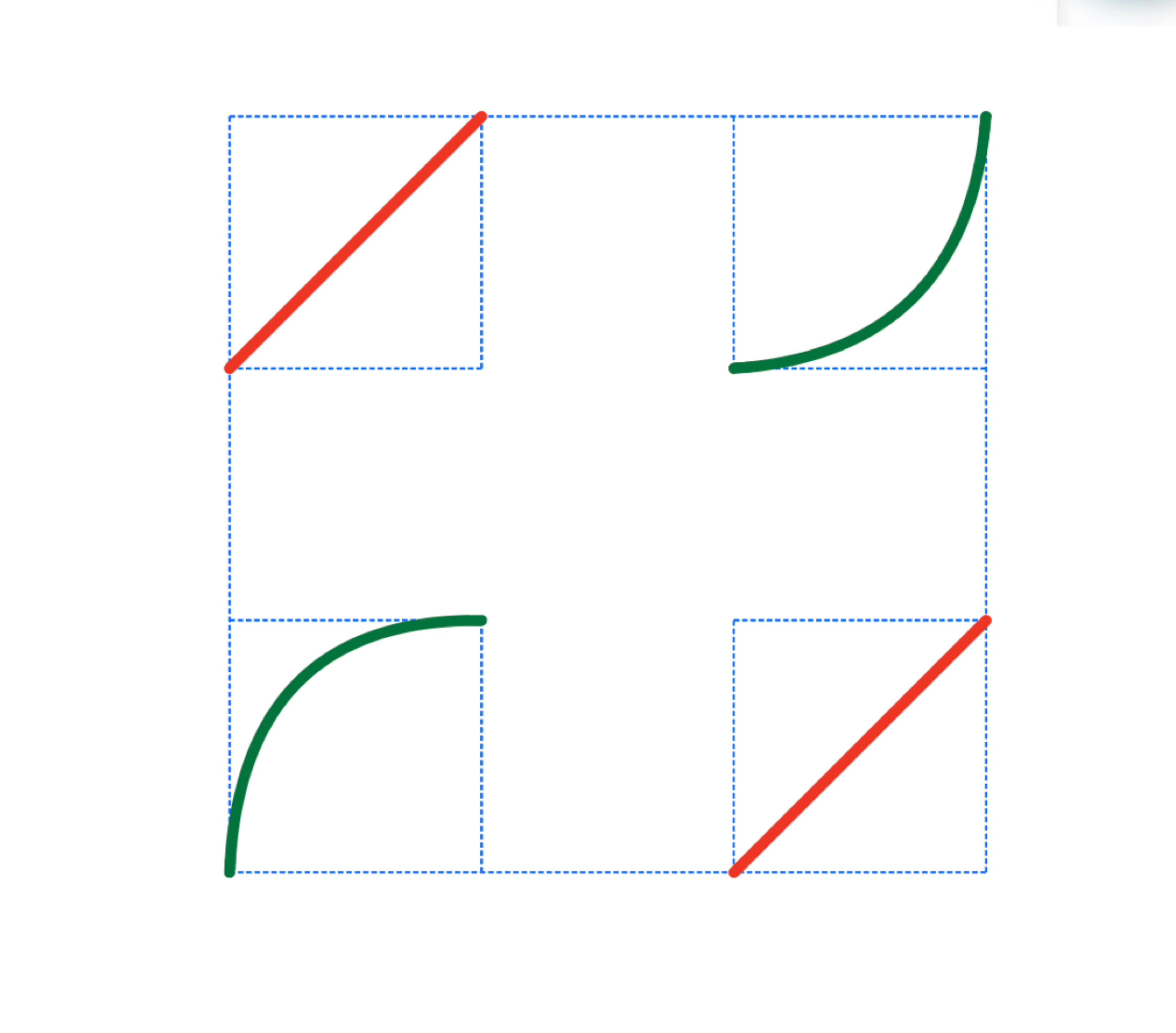}
	\caption{The relation $F$ from Example \ref{exx2}}
	\label{fig2}
\end{figure}  
 
Let $D=(0,1)\cup (2,3)$.  Note that $D$ is dense in $X$ and that for each $s\in D$,  $\Cl(\mathcal U^{\oplus}_F(s))=X$ (a slight modification of Lemma \ref{lemissimamissima} would give this result - we leave the details for the reader).  Therefore,  by Theorem \ref{main21}, $\sigma_F^+$ is transitive. It follows from Theorem \ref{main22} that $X_F$ is homeomorphic to $X_F^+$, and from Theorem \ref{tazadnji}  that $\sigma_F$ is a transitive homeomorphism on $X_F$.  

Next, we use this to obtain a transitive homeomorphism on the Cantor fan. 
For each $\alpha:\mathbb Z\rightarrow \{1,2\}$ and for each $t\in X$, let 
$
g(\alpha,t)=(\ldots,t_{-2},t_{-1} {,}t_0 {;} t_1,t_2,\ldots),
$ 
where $t_0=t$ and for each integer $k$, $t_k=f_{\alpha(k)}(t_{k-1})$. Note that $g$ is a homeomorphism from $C\times ([0,1]\cup [2,3])$ onto $X_F$, where $C=\{1,2\}^{\mathbb Z}$ is a Cantor set. For all $\mathbf x,\mathbf y\in X_F$, we define 
$$
\mathbf x\sim \mathbf y ~~~   \Longleftrightarrow  ~~~  \mathbf x = \mathbf y \textup{ or  for each integer } k,  {\{\mathbf x(k),\mathbf y(k)\}\subseteq} \{0,2\}. 
$$
Note that $X_F{/}_{\sim}$ is a Cantor fan and that for all $\mathbf x,\mathbf y\in X_F$, $\mathbf x\sim \mathbf y $ if and only if $\sigma_F(\mathbf x)\sim \sigma_F(\mathbf y)$. By Theorem \ref{kvocienti}, the function $\sigma_F^{\star}:X_F{/}_{\sim}\rightarrow X_F{/}_{\sim}$,  defined by 
$
\sigma_F^{\star}([x])=[\sigma_F(x)]
$
for each $x\in X$, is a transitive homeomorphism. 
}
\end{example}
\begin{example}\label{exx3}
{  (Our fourth proof of Theorem \ref{main})

Let $X=[-1,1]$ and let $F=\Gamma(f_1)\cup \Gamma(f_2)$, where $f_1:[-1,1]\rightarrow [-1,1]$ and $f_2:[-1,1]\rightarrow [-1,1]$ are defined by 
$$
f_1(x)=-x    ~~~  \textup{ and }  ~~~  f_2(x)=\begin{cases}
				x^{\frac{1}{3}}\text{;} & x\in [-1,0]\\
				x^2\text{;} & x\in [0,1]
			\end{cases} 
$$
for each $x\in X$; see Figure \ref{fig3}.
\begin{figure}[h!]
	\centering
		\includegraphics[width=15em]{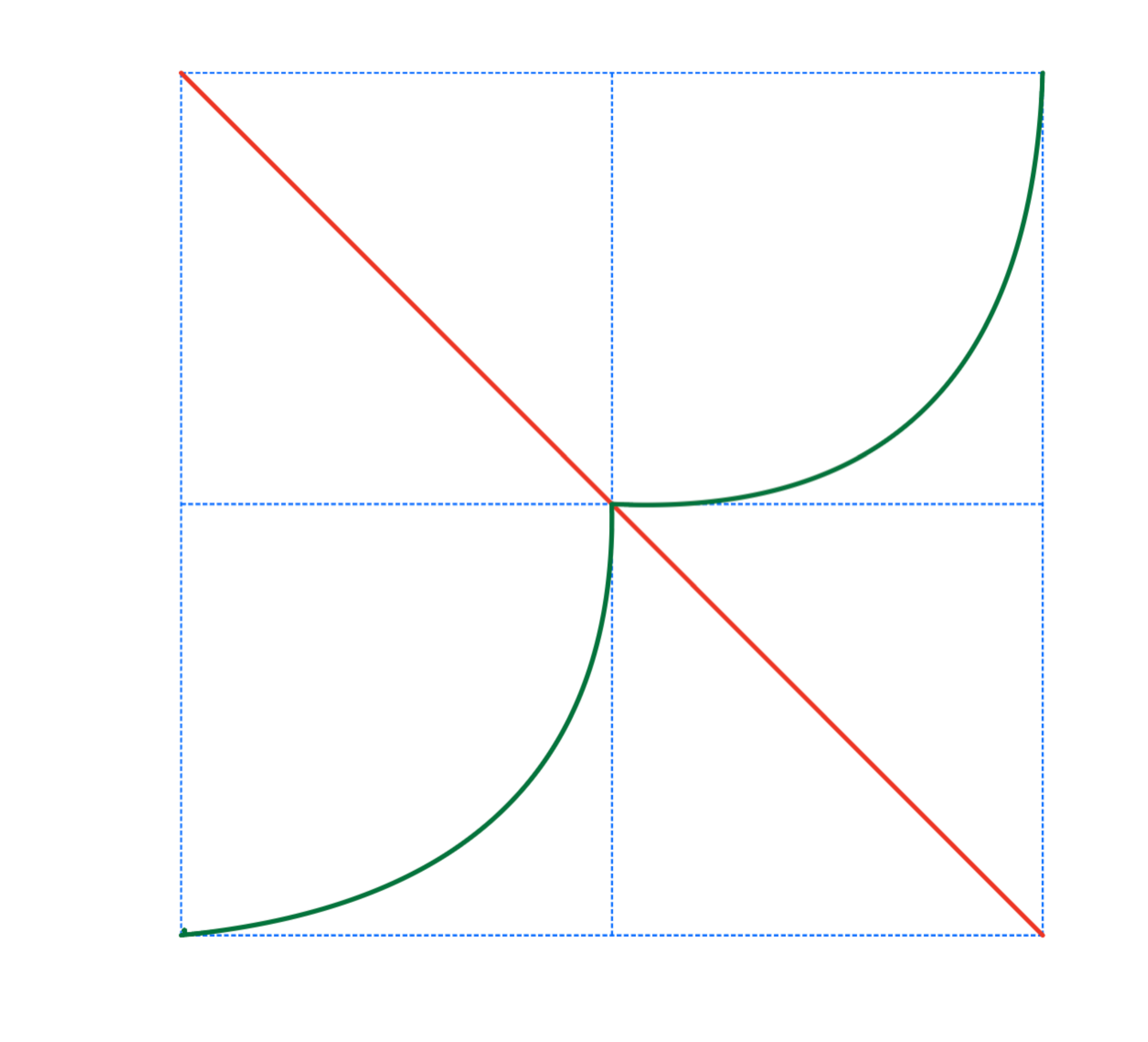}
	\caption{The relation $F$ from Example \ref{exx3}}
	\label{fig3}
\end{figure}  
 
 Let $D=(-1,0)\cup (0,1)$. Note that $D$ is dense in $X$ and that for each $s\in D$,   $\Cl(\mathcal U^{\oplus}_F(s))=X$ (a slight modification of Lemma \ref{lemissimamissima} would give this result - we leave the details for the reader).  Therefore,  by Theorem \ref{main21}, $\sigma_F^+$ is transitive. It follows from Theorem \ref{main22} that $X_F$ is homeomorphic to $X_F^+$, and from Theorem \ref{tazadnji}  that $\sigma_F$ is a transitive homeomorphism on $X_F$. Finally, we show that $X_F$ is a Cantor fan by showing that $X_F^+$ is a Cantor fan.  By \cite[Example 2.7, page 22]{ingram}, $X_G^+$, where $G=\{(x,x) \ | \ x\in [-1,1]\}\cup \{(x,-x) \ | \ x\in [-1,1]\}$, is a Cantor fan. Obviously, $X_F^+$ and $X_G^+$ are homeomorphic. Therefore, $X_F^+$ is a Cantor fan.
}
\end{example}

\section{Sigma transitivity}\label{s4}
In this section, we state and prove our second main result - Theorem \ref{mainn}. First, we introduce the notion of $\sigma$-transitivity of CR-dynamical systems and show that under certain minor conditions,  the CR-dynamical system $(X,F)$ is $\sigma$-transitive if and only if the dynamical system $(X_F^+,\sigma_F^+)$ is transitive. First, we define a $\sigma$-transitive CR-dynamical system.
\begin{definition}
Let $X$ be a compact metric space and let $F$ be a closed relation on $X$. Then we say that $(X,F)$ is \emph{a dynamical system with a closed relation} or \emph{ a CR-dynamical system}. 
\end{definition}
\begin{definition}\label{def}
Let $(X,F)$ be a CR-dynamical system. We say that $(X,F)$ is \emph{ $\sigma$-transitive}, if  there is a countable collection $\mathcal U=\{U_1,U_2,U_3,\ldots\}$
 of open sets in $\prod_{i=1}^{\infty}X$ such that
\begin{enumerate}
\item\label{1} for each positive integer $j$,  there is a positive integer $m_j$ and there are open sets $U_1^j$, $U_2^j$, $U_3^j$, $\ldots$, $U_{m_j}^j$  in $X$ such that
$
U_j=U_1^j\times U_2^j\times U_3^j\times  \ldots \times U_{m_j}^j \times \prod_{i=m_j+1}^{\infty}X,
$
\item\label{2} for each open set $V$ in $\prod_{i=1}^{\infty}X$ such that $V\cap X_F^+ \neq \emptyset$ and for each $\mathbf t\in V\cap X_F^+$, there is a positive integer $j$ such that 
$\mathbf t\in U_j\subseteq V$,
\item\label{3}  there is a sequence $(\mathbf x_j)$  of points in $X_F^+$ such that for each positive integer $j$,
\begin{enumerate}
\item $\mathbf x_j\in U_j$, and
 
\item  there is a positive integer $\ell_j > m_j$ such that
$
\pi_{\ell_j}(\mathbf x_j)=\pi_{1}(\mathbf x_{j+1}).
$
\end{enumerate}
\end{enumerate}
\end{definition}
 {We use the following well-known observation to prove Lemma \ref{lemma1}
 \begin{observation}
 Let $X$ be a compact metric space with cardinality greater than 1. Then the topological product $\Pi_{i=1}^\infty X$ does not contain any isolated points.
 \end{observation}}
\begin{lemma}\label{lemma1}
Let $(X,F)$ be a CR-dynamical system with the cardinality of $X$ greater than 1.  If $(X,F)$ is $\sigma$-transitive with $\mathcal U=\{U_1,U_2,U_3,\ldots\}$,  as in  Definition \ref{def},  then for each positive integer $j$,  there are infinitely many positive integers $k$ such that $U_k \subset U_j$.

\end{lemma}

\begin{proof}
Assume $(X,F)$ is a CR-dynamical system with the cardinality of  $X$ greater than 1, and assume $(X,F)$ is $\sigma$-transitive with $\mathcal U=\{U_1,U_2,U_3,\ldots\}$,  as in  Definition \ref{def} .  Suppose $j$ is a positive integer and $\mathbf t \in U_j \cap X_F^+ $ . Since no element of $\Pi_{i=1}^\infty X$ is isolated in $\Pi_{i=1}^\infty X$, there is an element $\mathbf y_1 \in U_j \setminus \{ \mathbf t \}$ . Now, $U_j \setminus \{ \mathbf y_1 \}$ is an open set such that $\mathbf t \in U_j \setminus \{ \mathbf y_1 \}$, thus there is a positive integer $j_1 \neq j$ such that $\mathbf t \in U_{j_1} \subset U_j $.  Similarly $U_j \setminus \{ \mathbf y_1 , \mathbf t \}$ is a non-empty open set in $\Pi_{i=1}^\infty X$. Let $ \mathbf y_2 \in U_j \setminus \{\mathbf y_1 , \mathbf t \} $. Then there is a positive integer $j_2 \not \in \{j, j_1 \} $ such that $\mathbf t \in U_{j_2} \subset U_{j_1} \setminus \{ \mathbf y_1,  \mathbf y_2 \} \subset U_j $.  A simple induction argument finishes the proof.
\end{proof}

 \begin{definition}
Let $X$ be a compact metric space and let $F$ be a closed relation on $X$.  Also, let $m$ and $n$ be positive integers and let $(x_{ 1}, x_{ 2},x_3,\ldots ,x_{ n}) \in \prod_{k=1}^{n}X$ and $(y_{ 1},y_{ 2}, y_3,\ldots , y_{ m} ) \in \prod_{k=1}^{m}X$ be such that  $x_{ n}=y_{ 1}$.  Then  we define \emph{ $(x_{ 1}, x_{ 2},x_3,\ldots , x_{ n})\star (y_{ 1}, y_{ 2}, y_{ 3}, \ldots , y_{ m})$} by 
$$
(x_{ 1}, x_{ 2}, x_{3},\ldots , x_{ n})\star (y_{1}, y_{2},y_{3}, \ldots , y_{ m}) = (x_{1}, x_{ 2},x_{ 3},\ldots , x_{ n}, y_{ 2}, y_{ 3},\ldots , y_{ m}).
$$
\end{definition}
The following observation easily follows.
\begin{observation}\label{lemma2}
Let $(X,F)$ be a $\sigma$-transitive  CR-dynamical system with $\mathcal U=\{U_1,U_2,U_3,\ldots\}$, $(\mathbf x_j)$ , $(m_j)$, and $(\ell_j)$ as in  Definition \ref{def}.   Also,  let  $s_0=0$ and  $s_k= \Sigma_{j=1}^k \ell _j -1 $ for {each positive integer $k$}, and let
$$
\mathbf x=\pi_{[1,\ell_1]}(\mathbf x_1)\star  \pi_{[1,\ell_2]}(\mathbf x_2)\star \pi_{[1,\ell_3]}(\mathbf x_3)\star  \ldots 
$$
Then $\mathbf x \in X_F^+$   and  $(\sigma_F^+)^{s_k}(\mathbf x) \in U_k$ for each non-negative integer $k$.
\end{observation}

We also need the well-known result from Observation \ref{jurc}.
\begin{observation}\label{jurc}
Let $(X,f)$ be a dynamical system. The following statements are equivalent.
\begin{enumerate}
\item \label{mili1} $(X,f)$ is  transitive.
\item \label{mili3} There is a point $x\in X$ such that $f(x)$ is a transitive point in $(X,f)$.
\end{enumerate}
\end{observation}

\begin{theorem}\label{main1}
Let $(X,F)$ be a CR-dynamical system.   If $(X,F)$ is $\sigma$-transitive, then $(X_F^+,\sigma_F^+)$ is transitive.
\end{theorem}

\begin{proof}
Let $(X,F)$ be a CR-dynamical system such that $(X,F)$ is $\sigma$-transitive with $\mathcal U=\{U_1,U_2,U_3,\ldots\}$, $(\mathbf x_j)$ , $(m_j)$, and $(\ell_j)$ as in  Definition \ref{def}. Also let  $s_0=0$ and  $s_k= \Sigma_{j=1}^k \ell _j -1 $ for each positive integer $k$, and let
$$
\mathbf x=\pi_{[1,\ell_1]}(\mathbf x_1)\star  \pi_{[1,\ell_2]}(\mathbf x_2)\star \pi_{[1,\ell_3]}(\mathbf x_3)\star  \ldots , 
$$
Also, let $\sigma=\sigma_F^+$. 
By {Observation} \ref{jurc} it suffices to show that  $\sigma^{\ell_1 -1}(\mathbf x)$ is a transitive point in $(X_F^+ , \sigma )$. To that end, let $V$ be an open set in $\Pi_{i=1}^\infty X$ such that $V \cap X_F^+ \neq \emptyset $. Then there is a positive integer $n$ such that $ U_n \subset V$.  By Lemma \ref{lemma1} there is a positive integer $k$ such that $k>n$ and $U_k \subset U_n$. Therefore, by {Observation} \ref{lemma2}, $\sigma^{s_k}(\mathbf x) \in U_k$.  Now $s_k - (\ell_1-1) > 0$ and $\sigma^{s_k}(\mathbf x) =\sigma^{s_k - (\ell_1-1)}(\sigma^{\ell_1 -1}(\mathbf x)) \in U_k \subset U_n$.  Thus we have shown that $\sigma^{\ell_1 -1}(\mathbf x)$ is a transitive point in $(X_F^+ , \sigma)$.
\end{proof}

Next, we present an example showing  that 
there is a continuum $X$ and  a closed relation $F$ on $X$ such that
\begin{enumerate}
\item $X_F^+$ has an isolated point,
\item $(X_F^+,\sigma_F^+)$ is dense orbit transitive, and
\item $(X,F)$ is not $\sigma$-transitive.
\end{enumerate}

\begin{example}
Let $X=[0,1]$, let $F=\{(0,1),(1,1)\}$ and let $\sigma=\sigma_F^+$.   Then 
$
\star_{i=1}^{\infty}F=\{(0,1,1,1,\ldots),(1,1,1,1,\ldots)\}
$
 and $\sigma(0,1,1,1,\ldots)=\sigma(1,1,1,1,\ldots)=(1,1,1,\ldots)$. The point  $(0,1,1,1,\ldots)$ is  transitive  in $(\star_{i=1}^{\infty}F,\sigma)$ since its orbit equals $\star_{i=1}^{\infty}F$:
 $$
 \{(0,1,1,1,\ldots),\sigma(0,1,1,1,\ldots),\sigma^2(0,1,1,1,\ldots),\ldots\}=\{(0,1,1,1,\ldots),(1,1,1,1,\ldots)\}
 $$
 Next, we show that $(X,F)$ is not $\sigma$-transitive. Suppose $(X,F)$ is $\sigma$-transitive with $\mathcal U=\{U_1,U_2,U_3,\ldots\}$  and sequence $(m_j)$ as in  Definition \ref{def}. We show that $(X,F)$ is not $\sigma$-transitive by showing that {\bf it does not hold} that there are a sequence $(\mathbf x_j)$ of points in $X_F^+$ and a sequence of positive integers $(\ell_j)$ such that for each positive integer $j$, $\ell_j >m_j$,  $\mathbf x_j\in U_j$, and 
$\pi_{\ell_j}(\mathbf x_j)=\pi_1(\mathbf x_j)$.

 Let 
$$
\mathcal U_0=\{U\in \mathcal U \ | \ (0,1,1,1,\ldots)\in U, (1,1,1,1,\ldots)\not \in U\}
$$
and let 
$$
\mathcal U_1=\{U\in \mathcal U \ | \ (0,1,1,1,\ldots)\not \in U, (1,1,1,1,\ldots) \in U\}. 
$$
Note that each of the {collections} $\mathcal U_0$ and $\mathcal U_1$ is infinite. Therefore, there is a positive integer $j$ such that $U_j\in \mathcal U_1$ and $U_{j+1}\in \mathcal U_0$. It follows that  if  $\mathbf x_j \in U_j\cap X_F^+$, then $\mathbf x_j =(1,1,1,1,\ldots)$ and, therefore, for each positive integer $\ell$, there is no $\mathbf x_{j+1} \in U_{j+1} \cap X_F^+ $ such that $\pi_1(\mathbf x_{j+1})=\pi_{\ell}(\mathbf x_j)$  since $\pi_{\ell}(\mathbf x_j)=1$ while $\pi_1(\mathbf x_{j+1})=0$.
\end{example}

We use the well-known result from Observation \ref{lemissima1} to prove the converse of Theorem \ref{main1}.

\begin{observation}\label{lemissima1}
Let $(X,f)$ be a dynamical system.  If $(X,f)$ is transitive, then for each transitive point $x$ in $(X,f)$ and for each positive integer $n$, $f^n(x)$ is also a transitive point in $(X,f)$. 
\end{observation}

\begin{theorem}\label{prva}
Let $(X,F)$ be a CR-dynamical system.  If $(X_F^+,\sigma_F^+)$ is transitive, then $(X,F)$ is $\sigma$-transitive.
\end{theorem}
\begin{proof}
Let $\sigma=\sigma_F^+$ and suppose that $(X_F^+,\sigma)$ is transitive.  By Observation \ref{isolatedpoints},  the dynamical system $(X_F^+,\sigma)$ is dense orbit transitive. Let
$
\mathbf x=(x_1,x_2,x_3,\ldots)\in X_F^+
$
 be a transitive point in $(X_F^+,\sigma)$ and  let $\mathcal B$ be a countable base for $\prod_{i=1}^{\infty}X$ such that for each $B\in \mathcal{B}$, there are a positive integer $m$ and open sets $B_1$, $B_2$, $B_3$, $\ldots$, $B_{m}$  in $X$ such that 
$
B=B_1\times B_2\times B_3\times \ldots \times B_{m}\times \prod_{i=m+1}^{\infty}X. 
$
Also, let  
$$
\mathcal U=\{B\in \mathcal B \ | \ B \cap \star_{i=1}^{\infty}F\neq\emptyset\}=\{U_1,U_2,U_3,\ldots\}. 
$$
For each positive integer $j$, let $m_j$ be a positive integer and let $U_1^j$, $U_2^j$, $U_3^j$, $\ldots$,  $U_{m_j}^j$  be open sets in $X$ such that
$
U_j=U_1^j\times U_2^j\times U_3^j\times \ldots \times U_{m_j}^j \times \prod_{i=m_j+1}^{\infty}X.
$
Then the following hold.
\begin{enumerate}
\item Obviously, \ref{1} from Definition \ref{def} holds.
\item Let $V$ be any open set in $\prod_{i=1}^{\infty}X$ and let $\mathbf t\in V\cap X_F^+$.  Since $\mathcal B$ is a base for $\prod_{i=1}^{\infty}X$,  there is a positive integer $j$ such that $\mathbf t\in U_j\subseteq V$.  Therefore, \ref{2} from Definition \ref{def} is satisfied. 
\item We construct sequences $(\ell_j)$ and $(\mathbf x_j)$  satisfying \ref{3} from Definition \ref{def}  inductively.  
\begin{enumerate}
\item Let $j=1$.  Then let $k$ be such a positive integer that $\sigma^{k}(\mathbf x)\in U_1$.  Such a positive integer $k$ does exist since  $\mathbf x$ is  a transitive point in $(X_F^+,\sigma)$.  Let  
$
\mathbf x_1=\sigma_F^{k}(\mathbf x).
$

\item  Let $j=2$.  Then let $\ell_1$ be a positive integer such that $\ell_1 >m_1$ and $\sigma^{\ell_1-1}(\mathbf x_1)\in U_2$. Such a positive integer $\ell_1$ does exist since,  by Observation \ref{lemissima1}, $\mathbf \sigma^{m_1}(x_1)$ is also a transitive point in $(X_F^+,\sigma)$. Let 
$
\mathbf x_2=\sigma^{\ell_1-1}(\mathbf x_1). 
$
Note that $\ell_1 >m_1$ , $\pi_{\ell_1}(\mathbf x_1)=\pi_1(\mathbf x_2)$, $\mathbf x_1 \in U_1$, and $\mathbf x_2 \in U_2$.

\item Let $j>1$ be a positive integer and suppose that we have already constructed  positive integers $\ell_1 , \ell_2, \ldots , \ell_{j-1}$ and points $\mathbf x_{1}, \mathbf x_2, \ldots \mathbf x_{j}$ in $X_F^+$ so that for each $k < j $ we have $\ell_{k} >m_{k}$ , $\pi_{\ell_{k}}(\mathbf x_{k})=\pi_1(\mathbf x_{k+1})$, and $\mathbf x_k \in U_k$, and that we also have $\mathbf x_j \in U_j$. Then we construct $\ell_{j}$, and $\mathbf x_{j+1}$  as follows.   Let $\ell_{j}$ be a positive integer such that $\ell_j >m_j$ and $\sigma^{\ell_{j}-1}(\mathbf x_{j-1})\in U_{j}$.  Such a positive integer $\ell_j$ does exist since,  by Observation \ref{lemissima1}, $ \sigma^{m_{j-1}}(\mathbf x_{j-1})$ is also a transitive point in $(X_F^+,\sigma)$.  Let 
$
\mathbf x_{j+1}=\sigma^{\ell_j -1}(\mathbf x_{j}). 
$ 
Note that $\ell_j >m_j$ , $\pi_{\ell_j}(\mathbf x_j)=\pi_1(\mathbf x_{j+1})$ , and $\mathbf x_{j+1} \in U_{j+1}$.
\end{enumerate}
The sequences $(\ell_j )$ and $(\mathbf x_j)$   were constructed in such a way that  \ref{3} from Definition \ref{def} is satisfied. 
\end{enumerate}
It follows that $(X,F)$ is $\sigma$-transitive.   
\end{proof}
Our next goal is to prove Theorem \ref{mainn}, which is one of our main results. First, we prove the following lemma that will be needed in its proof.
\begin{lemma}\label{Gabi}
For each $x\in (0,1]$, the set 
$$
\left\{\Big(\frac{1}{2}\Big)^{\frac{3^{h}-1}{2^{k+1}}}\cdot x^{\frac{3^{h}}{2^k}} \ \Big| \  k \textup{ and } h \textup{ are positive integers}\right\}
$$
is dense in $[0,1]$. 
\end{lemma}
\begin{proof}
We prove the lemma by showing that 
$$
\left\{\ln\left(\Big(\frac{1}{2}\Big)^{\frac{3^{h}-1}{2^{k+1}}}\cdot x^{\frac{3^{h}}{2^k}}\right) \ \Big| \  k \textup{ and } h \textup{ are positive integers}\right\}
$$
is dense in $(-\infty,0]$.  First, note that
$$ 
\ln\left(\Big(\frac{1}{2}\Big)^{\frac{3^{h}-1}{2^{k+1}}}\cdot x^{\frac{3^{h}}{2^k}}\right)=\ln\Big(\frac{1}{2}\Big)^{\frac{3^{h}-1}{2^{k+1}}}+\ln x^{\frac{3^{h}}{2^k}}={\frac{3^{h}-1}{2^{k+1}}}\cdot \ln\frac{1}{2}+{\frac{3^{h}}{2^k}}\cdot \ln x
$$
and let $a=\ln\frac{1}{2}$ and $b=\ln x$. Then $a<0$ and $b<0$, and it follows that 
$$
\ln\left(\Big(\frac{1}{2}\Big)^{\frac{3^{h}-1}{2^{k+1}}}\cdot x^{\frac{3^{h}}{2^k}}\right)={\frac{3^{h}-1}{2^{k+1}}}\cdot a+{\frac{3^{h}}{2^k}}\cdot b={\frac{3^{h}}{2^{k}}}\cdot \frac{a}{2}+{\frac{3^{h}}{2^k}}\cdot b-\frac{a}{2^{k+1}}={\frac{3^{h}}{2^{k}}}\cdot \Big(\frac{a}{2}+b\Big)-\frac{a}{2^{k+1}}.
$$
To conclude the proof, we show that for all $a<0$ and $b<0$,
$$
\left\{{\frac{3^{h}}{2^{k}}}\cdot \Big(\frac{a}{2}+b\Big)-\frac{a}{2^{k+1}} \ \Big| \  k \textup{ and } h \textup{ are positive integers}\right\}
$$
is dense in $(-\infty,0]$.  Let $c<d\leq 0$ and let $d'=\frac{c+d}{2}$.   A similar argument as in the proof of Lemma \ref{lemissimamissima} proves that the set  $\big\{\frac{3^h}{2^k} \ | \  h \textup{ and } k \textup{ are positive integers}\big\}$ is dense in $[0,\infty)$. 
Let $k$ and $h$ be such positive integers that 
$$
\frac{d'}{\frac{a}{2}+b}<\frac{3^h}{2^k}<\frac{c}{\frac{a}{2}+b}  \textup{ and } ~~~  \frac{-a}{2^{k+1}}<\frac{d-c}{2}.
$$ 
Then $c<\frac{3^h}{2^k}(\frac{a}{2}+b)<d'$, and it follows that
$$
c<c-\frac{a}{2^{k+1}}<\frac{3^h}{2^k}(\frac{a}{2}+b)-\frac{a}{2^{k+1}}<d'-\frac{a}{2^{k+1}}{<}\frac{c+d}{2}+\frac{d-c}{2}=d.
$$
This completes the proof.
\end{proof}

\begin{definition}
We  use  $I$ to denote $I=[0,1]$ and $H$ to denote 
$$
H=\big\{\big(x,\sqrt x\big) \ | \ x\in [0,1]\big\}\cup \Big\{\Big(x,\frac{1}{2}x^3\Big) \ | \ x\in [0,1]\Big\}.
$$
We also use $f_0$ and $f_1$ to denote the functions $f_0,f_1:I\rightarrow I$ that are defined by $f_0(x)=\frac{1}{2}x^3$ and $f_1(x)=\sqrt x$ for any $x\in I$. 
\end{definition}
\begin{theorem}\label{Cantor fan plus}
$I_H^+$ is a Cantor fan.
\end{theorem}
\begin{proof}
It follows from \cite[Example 1, page 7]{banic1} that $I_G^+$, where $G=\{(x,x) \ | \ x\in I\}\cup \{(x,\frac{1}{2}x) \ | \ x\in I\}$ is a Cantor fan.  Obviously, the function $\varphi:I_G^+\rightarrow I_H^+$, defined by 
$$
\varphi(x_1,x_2,x_3,\ldots)=(x_1,g_1(x_1),g_2(g_1(x_1)),g_3(g_2(g_1(x_1))),\ldots)
$$
for any $(x_1,x_2,x_3,\ldots)\in I_G^+$,  where for each positive integer $k$, $g_k=f_0$ if $x_{k+1}=\frac{1}{2}x_k$ and $g_k=f_1$ if $x_{k+1}=x_k$, is a homeomorphism. 
\end{proof}
\begin{observation}\label{JudyGoran}
For each $x\in (0,1]$,
$
f_1^k(f_0^h(x))=\Big(\frac{1}{2}\Big)^{\frac{3^{h}-1}{2^{k+1}}}\cdot x^{\frac{3^{h}}{2^k}}.
$
\end{observation}
\begin{theorem}\label{jojmene1}
The shift map $\sigma_{H}^{+}$  is a transitive continuous surjection
\end{theorem}
\begin{proof}
$\sigma_{H}^{+}$  is a  continuous surjection. To show that it is  transitive, we show that the CR-dynamical system $(I,H)$ {is} $\sigma$-transitive.   To do this, let $\mathcal B$ be a countable base for $\prod_{i=1}^{\infty}I$ such that for each $B\in \mathcal{B}$, there are a positive integer $m$ and open sets $B_1$, $B_2$, $B_3$,  $\ldots$, $B_{m}$ in $I$ such that 
$
B=B_1\times B_2\times B_3\times \ldots \times B_{m}\times \prod_{i=m+1}^{\infty}I. 
$
Also, let  
$$
\mathcal U=\{B\in \mathcal B \ | \ B \cap I_H^+\neq\emptyset\}=\{U_1,U_2,U_3,\ldots\}. 
$$
For each positive integer $j$, let $m_j$ be a positive integer and let $U_1^j$, $U_2^j$, $U_3^j$, $\ldots$, $U_{m_j}^j$ be open sets in $I$ such that
$
U_j=U_1^j\times U_2^j\times U_3^j\times \ldots \times U_{m_j}^j\times \prod_{i=m_j+1}^{\infty}I.
$
Then the following hold.
\begin{enumerate}
\item  Obviously,  \ref{1} from Definition \ref{def} is satisfied. 
\item Let $V$ be any open set in $\prod_{i=1}^{\infty}I$ and let $\mathbf t\in V\cap I_H^+$.  Since $\mathcal B$ is a base for $\prod_{i=1}^{\infty}I$ and since $\mathbf t\in I_H^+$,  there is a positive integer $j$ such that $\mathbf t\in U_j\subseteq V$.
Therefore, \ref{2} from Definition \ref{def} is satisfied. 
\item Let $j$ be a positive integer and let 
$
\mathbf x_j'=(x_1^j,x_2^j,x_3^j,\ldots)\in \Big(U_j\cap I_H^+\Big)\setminus\{(0,0,0,\ldots)\}.
$
Also, let $(f_1^j,f_2^j,f_3^j,\ldots)\in \{f_0,f_1\}^{\mathbb N}$ be such a sequence that for each positive integer $n$,
$
x_{n+1}^j=f_n^j(x_{n}^j).
$
Next,  let  $V_{m_j}^{j-1}$ be an open set in $[0,1]$ such that 
$
x_{m_j-1}^j\in V_{m_j-1}^{j}\subseteq U_{m_j-1}^{j}
$
 and for each $x\in V_{m_j-1}^j$, $f_{m_j-1}^{j}(x)\in U_{m_j}^{j}$. Note that such a set  $V_{m_j-1}^{j}$ does exist since $f_{m_j-1}^j$ is continuous.  Next, let $V_{m_j-2}^{j}$ be an open set in $I$ such that 
 $
 x_{m_j-1}^j\in V_{m_j-1}^{j}\subseteq U_{m_j-1}^{j}
 $
  and for each $x\in V_{m_j-1}^{j}$, $f_{m_j-1}^j(x)\in V_{m_j}^{j}$. Again, such a set  $V_{m_j-1}^{j}$ does exist since $f_{m_j-1}^j$ is continuous.  We continue inductively. Finally, suppose that we have already constructed the sets $V_{m_j}^{j}$,  $V_{m_j-1}^{j}$,  $V_{m_j-2}^{j}$,  $\ldots$, $V_{2}^{j}$.  Then we construct the set $V_{1}^{j}$ as follows.  
Let $V_{1}^{j}$ be an open set in $[0,1]$ such that 
$
x_1^j\in V_{1}^{j}\subseteq U_{1}^{j}
$
 and for each $x\in V_{1}^{j}$, $f_1^j(x)\in V_{2}^{j}$. Such a set  $V_{1}^{j}$ does exist since $f_{1}^j$ is continuous. 
 
Next, we define the sequence $(\mathbf x_k)$ of points in $I_H^+$ such that \ref{3} from Definition \ref{def} is satisfied as follows.  
\begin{itemize}
\item Let $x_1\in V_1^1\cap (0,1]$ and let $y_1=f_{m_1-1}^1(\ldots f_3^1(f_2^1(f_1^1(x_{{1}})))\ldots)$. Then $y_1\in V_{m_1}^1{\cap (0,1]}$. By Lemma \ref{Gabi}, the set 
$$
\left\{\Big(\frac{1}{2}\Big)^{\frac{3^{h}-1}{2^{k+1}}}\cdot y_1^{\frac{3^{h}}{2^k}} \ \Big| \  k \textup{ and } h \textup{ are positive integers}\right\}
$$
is dense in $[0,1]$.  Therefore, there are positive integers $k_1$ and $h_1$  such that 
$
\Big(\frac{1}{2}\Big)^{\frac{3^{h_1}-1}{2^{k_1+1}}}\cdot y_1^{\frac{3^{h_1}}{2^{k_1}}}\in V_1^2.
$
 Fix such positive integers $k_1$ and $h_1$. Then let 
	\begin{align*}
	 \mathbf y_1=\Big(&x_1, f_1^1(x_1), f_2^1(f_1^1(x_1)),\ldots , f_{m_1-{2}}^1(\ldots f_3^1(f_2^1(f_1^1(x_1)))\ldots){,}y_1,
	 \\
	 &f_0(y_1),f_0^2(y_1),f_0^3(y_1),\ldots,f_0^{h_1}(y_1),
	 \\
	 &f_1(f_0^{h_1}(y_1)),f_1^2(f_0^{h_1}(y_1)),f_1^3(f_0^{h_1}(y_1)),\ldots ,f_1^{k_1}(f_0^{h_1}(y_1))\Big),
	\end{align*}
	let $\ell_1=m_1+h_1+k_1$, and let  
	$\mathbf x_1\in I_H^+$ be any point such that 
	$
	\pi_{[1,\ell_1]}(\mathbf x_1)=\mathbf y_1.
	$
	{Note that $\mathbf y_1$ is an element of the finite Mahavier product $ I_H^{m_1+h_1+k_1-1}$. } Also, note that $\mathbf x_1\in U_1$ and that $\pi_{\ell_1}(\mathbf x_1)\in V_1^2\cap(0,1]$.
	\item Let $x_2=\pi_{\ell_1}(\mathbf x_1)$ and let $y_2=f_{m_2-1}^2(\ldots f_3^2(f_2^2(f_1^2(x_{2})))\ldots)$. Then $y_2\in V_{m_2}^2{\cap (0,1]}$. By Lemma \ref{Gabi}, the set 
$$
\left\{\Big(\frac{1}{2}\Big)^{\frac{3^{h}-1}{2^{k+1}}}\cdot y_2^{\frac{3^{h}}{2^k}} \ \Big| \  k \textup{ and } h \textup{ are positive integers}\right\}
$$
is dense in {$I$}.  Therefore, there are positive integers $k_2$ and $h_2$  such that 
$
\Big(\frac{1}{2}\Big)^{\frac{3^{h_1}-1}{2^{k_1+1}}}\cdot y_2^{\frac{3^{h_1}}{2^{k_1}}}\in V_1^3.
$
 Fix such positive integers $k_2$ and $h_2$. Then let 
	\begin{align*}
	 \mathbf y_2=\Big(&x_2, f_1^2(x_2), f_2^2(f_1^2(x_2)),\ldots , f_{m_2-{2}}^2(\ldots f_3^2(f_2^2(f_1^2(x_2)))\ldots){,}y_2,
	 \\
	 &f_0(y_2),f_0^2(y_2),f_0^3(y_2),\ldots,f_0^{h_2}(y_2),
	 \\
	 &f_1(f_0^{h_2}(y_2)),f_1^2(f_0^{h_2}(y_2)),f_1^3(f_0^{h_2}(y_2)),\ldots ,f_1^{k_2}(f_0^{h_2}(y_2))\Big)\in I_H^{m_2+h_2+k_2-1},
	\end{align*}
	 let $\ell_2=m_2+h_2+k_2$, and let  
	$\mathbf x_2\in I_H^+$ be any point such that 
	$
	\pi_{[1,\ell_2]}(\mathbf x_2)=\mathbf y_2.
	$
	Note that $\mathbf x_2\in U_2$ and that $\pi_{\ell_2}(\mathbf x_2)\in V_1^3\cap(0,1]$.
\end{itemize}
Continuing inductively, we construct  the sequence $(\mathbf x_k)$ such that \ref{3} from Definition \ref{def} is satisfied. 
\end{enumerate}
It follows  that $(X,F)$ is $\sigma$-transitive.   
\end{proof}
{Note that $\sigma_H^+$ from Theorem \ref{jojmene1} is not a homeomorphism}. Next, we show that the inverse limit $\varprojlim(I_H^+,\sigma_H^+)$ is homeomorphic to the Lelek fan. We need the following definition.
{
\begin{observation}
Observe that for any $\mathbf x\in I_H\setminus \{(\ldots,0,0,0;0,0,\ldots)\}$, there are unique sequences $(a_1,a_2,a_3,\ldots), (b_1,b_2,b_3,\ldots)\in \{f_{0},f_{1}\}^{\mathbb N}$ and a unique point $t\in (0,1]$ such that 
$$
\mathbf x=(\ldots ,a_3^{-1}(a_2^{-1}(a_1^{-1}(t))), a_2^{-1}(a_1^{-1}(t)),a_1^{-1}(t); t,b_1(t),b_2(b_1(t)),b_3(b_2(b_1(t))),\ldots).
$$
\end{observation}
\begin{definition}
Let $g_0$ and $g_1$ be the functions $g_0,g_1:[0,\infty)\rightarrow [0,\infty)$, defined by $g_0(x)=\frac{1}{2}x^3$ and $g_1(x)=\sqrt x$ for any $x\in [0,\infty)$. For any 
$
\mathbf h=(h_1,h_2,h_3,\ldots)\in \{g_{0},g_{1}\}^{\mathbb N},
$
we use 
$$
\overrightarrow{\mathbf h}(t)=(t,h_1(t),h_2(h_1(t)),h_3(h_2(h_1(t))),\ldots)
$$
and 
$$
\overleftarrow{\mathbf h}(t)=(\ldots ,h_{3}^{-1}(h_{2}^{-1}(h_{1}^{-1}(t))),h_{2}^{-1}(h_{1}^{-1}(t)),h_{1}^{-1}(t))
$$
for any $t\in I$. 
\end{definition}
\begin{observation}
Observe that for any $\mathbf x\in I_H\setminus \{(\ldots,0,0,0;0,0,\ldots)\}$, there are unique sequences $\mathbf a, \mathbf b\in \{g_{0},g_{1}\}^{\mathbb N}$ and a unique point $t\in (0,1]$ such that 
$$
\mathbf x=\overleftarrow{\mathbf a}(t)\oplus \overrightarrow{\mathbf b}(t).
$$
\end{observation}}
\begin{lemma}
$I_H$ is a subcontinuum of a Cantor fan.
\end{lemma}
\begin{proof}
We show that $I_H$ is a fan by showing that $I_H$ is a subcontinuum of a Cantor fan.  Let $g_0$ and $g_1$ be the functions $g_0,g_1:[0,\infty)\rightarrow [0,\infty)$, defined by $g_0(x)=\frac{1}{2}x^3$ and $g_1(x)=\sqrt x$ for any $x\in [0,\infty)$.  For all {$\mathbf a,\mathbf b\in \{g_0,g_1\}^{\mathbb N}$}, let 
$$
A_{\mathbf a,\mathbf b}=\{\overleftarrow{\mathbf a}(t)\oplus \overrightarrow{\mathbf b}(t) \ | \ t\in [0,1]\}. 
$$
{Note that for all $\mathbf a,\mathbf b\in \{g_0,g_1\}^{\mathbb N}$, $A_{\mathbf a,\mathbf b}$ is an arc  with end-points $(\ldots,0,0,0;0,0,\ldots)$ and $\mathbf e_{\mathbf a,\mathbf b}=\overleftarrow{\mathbf a}(1)\oplus \overrightarrow{\mathbf b}(1)$ in the Hilbert cube 
$$
\ldots \times \left[0,\sqrt[3]{2\sqrt[3]{2}}\right]\times \left[0,\sqrt[3]{2}\right]\times [0,1]\times [0,1]\times [0,1] \times \ldots.
$$}
Also, note that  for all $\mathbf a,\mathbf b,\mathbf a',\mathbf b'\in \{g_0,g_1\}^{\mathbb N}$, 
$$
(\mathbf a,\mathbf b)\neq (\mathbf a',\mathbf b') ~~~  \Longrightarrow ~~~   A_{\mathbf a,\mathbf b}\cap A_{\mathbf a',\mathbf b'}=\{(\ldots,0,0,0;0,0,\ldots)\}.
$$ 
Let $X=\bigcup_{(\mathbf a,\mathbf b)\in \{g_0,g_1\}^{\mathbb N}}A_{\mathbf a,\mathbf b}$ and let $E(X)=\{\mathbf e_{\mathbf a,\mathbf b} \ | \ (\mathbf a,\mathbf b)\in \{g_0,g_1\}^{\mathbb N}\}$.  Next, we show that $E(X)$ is homeomorphic to the Cantor set $\{g_0,g_1\}^{\mathbb N}\times \{g_0,g_1\}^{\mathbb N}$ (here, $\{g_0,g_1\}$ is equipped with the discrete topology while  $\{g_0,g_1\}^{\mathbb N}$ and $\{g_0,g_1\}^{\mathbb N}\times \{g_0,g_1\}^{\mathbb N}$  are equipped with the product topologies).   Let $\varphi:\{g_0,g_1\}^{\mathbb N}\times \{g_0,g_1\}^{\mathbb N}\rightarrow E(X)$ be defined by $\varphi(\mathbf a,\mathbf b)=\mathbf e_{\mathbf a,\mathbf b}$ for any $(\mathbf a,\mathbf b)\in \{g_0,g_1\}^{\mathbb N}\times \{g_0,g_1\}^{\mathbb N}$.  Note that $\varphi$ is a bijection. To see that $\varphi$ is a homeomorphism, we only need to see that it is continuous (since every continuous bijection from a compact space to a metric space is a homeomorphism).  To see that $\varphi$ is continuous, let $(\mathbf a,\mathbf b)\in \{g_0,g_1\}^{\mathbb N}\times \{g_0,g_1\}^{\mathbb N}$  and let $(\mathbf a_n,\mathbf b_n)$ be a sequence in $\{g_0,g_1\}^{\mathbb N}\times \{g_0,g_1\}^{\mathbb N}$ such that $\displaystyle \lim_{n\to \infty}(\mathbf a_n,\mathbf b_n)=(\mathbf a,\mathbf b)$.  Since the sequences $(\mathbf a_n)$ and $(\mathbf b_n)$ are coordinate-wise converging to $\mathbf a$ and $\mathbf b$,  respectively, it follows that $\displaystyle \lim_{n\to \infty}\mathbf e_{\mathbf a_n,\mathbf b_n}=\mathbf e_{\mathbf a,\mathbf b}$.  Therefore,  $\varphi $ is a homeomorphism and it follows that $X$ is a Cantor fan.  Since $I_H$ is a subcontinuum of $X$, it follows that $I_H$ is a subcontinuum of a Cantor fan.
\end{proof}
\begin{lemma}\label{totolele}
Let $\mathbf x\in I_H\setminus \{(\ldots,0,0,0;0,0,\ldots)\}$.  If there is an integer $k$ such that $\mathbf x(k)=1$, then $\mathbf x\in E(I_H)$.
\end{lemma}
\begin{proof}
Let $x\in (0,1]$,  let $\mathbf a, \mathbf b\in \{f_{0},f_{1}\}^{\mathbb N}$ such that $\mathbf x=\overleftarrow{\mathbf a}(x)\oplus \overrightarrow{\mathbf b}(x)$.   Suppose that $\mathbf x\not \in E(I_H)${. L}et $A$ be the maximal arc in $I_H$ with one end-point being $(\ldots,0,0,0;0,0,\ldots)$ such that $\mathbf x\in A$.  Also,  let $\pi_0(A)=[0,e]$.  It follows that $x\leq e$.  Since the point $\mathbf x$ is not an end-point of $I_H$, it follows that $x<e$.  Let $\mathbf e=\overleftarrow{\mathbf a}(e)\oplus \overrightarrow{\mathbf b}(e)$.  Note that  for each positive integer $k$,
$
\mathbf x(k)<\mathbf e(k).
$ 
If there is an integer $k$ such that $\mathbf x(k)=1$, then $\mathbf e(k)>1$, which is a contradiction. Therefore, for each integer $k$, $\mathbf x(k)<1$.  This completes the proof.
\end{proof}

\begin{theorem}\label{juuuj}
 The two-sided Mahavier product $I_H$ is a Lelek fan.
\end{theorem}
\begin{proof}
Let $\mathbf x=(\ldots,x_{-2},x_{-1};x_{0},x_{1},x_{2},\ldots)\in I_H$ be any point and let $\varepsilon >0$.  Also, let $m$ be a positive integer such that  $\sum_{k=m}^{\infty}\frac{1}{2^k}<\frac{\varepsilon}{3}$.  For each integer $k\in \{-m,-m+1,-m+2,\ldots, m-2,m-1,m\}$, let $g_k\in \{f_0,f_1\}$ be such that $x_{k+1}=g_k(x_k)$. Next, let for each integer $k\in \{-m,-m+1,-m+2,\ldots, m-2,m-1,m\}$, $V_k$ be an open set in $I$ such that 
$$
\mathbf x\in\Big(\prod_{k=-\infty}^{-m-1}I\Big)\times V_{-m}\times V_{-m+1}\times V_{-m+2}\times \ldots \times V_{m-2}\times V_{m-1}\times V_{m}\times \Big(\prod_{i=m+1}^{\infty}I\Big)\subseteq B(\mathbf x,\frac{\varepsilon}{3}).
$$ 
Without loss of generality, we may assume that for each integer $k\in \{-m,-m+1,-m+2,\ldots, m-2,m-1\}$,  $g_k(V_k)\subseteq V_{k+1}$.  

By Lemma \ref{Gabi},  the set $\big\{\big(\frac{1}{2}\big)^{\frac{3^{h}-1}{2^{k+1}}} \ \big| \  k \textup{ and } h \textup{ are positive integers}\big\}$ 
is dense in $[0,1]$ and by Observation \ref{JudyGoran},  $f_1^k(f_0^h(1))=\big(\frac{1}{2}\big)^{\frac{3^{h}-1}{2^{k+1}}}$.  
Let $h_0$ and $k_0$ be positive integers such that $f_1^{k_0}(f_0^{h_0}(1))\in V_{-m}$. We define $\mathbf e$ as follows:
\begin{enumerate}
\item $\mathbf e(-m-h_0-k_0)=1$,
\item for each  $k\leq -m-k_0-h_0$, let $\mathbf e(k-1)=f_1^{-1}(\mathbf e(k))$,
\item for each  $k\in \{-m-k_0-h_0,-m-k_0-h_0+1,-m-k_0-h_0+2,\ldots, -m-h_0-1\}$, let $\mathbf e(k+1)=f_0(\mathbf e(k))$,
\item for each  $k\in \{-m-h_0,-m-h_0+1,-m-h_0+2,\ldots, -m-1\}$, let $\mathbf e(k+1)=f_1(\mathbf e(k))$,
\item for each  $k\in \{-m,-m+1,-m+2,\ldots, m-2,m-1\}$, let $\mathbf e(k+1)=g_k(\mathbf e(k))$,
\item for each  $k\geq m$, let $\mathbf e(k+1)=f_0(\mathbf e(k))$.
\end{enumerate}
By Lemma  \ref{totolele}, $\mathbf e\in E(I_H)$ and it follows from its construction that 
$$
D_H(\mathbf x,\mathbf e)=\sum_{k=-\infty}^{-m-1}\frac{|x_k-\mathbf e(k)|}{2^{|k|}} +\sum_{k=-m}^{m}\frac{|x_k-\mathbf e(k)|}{2^{|k|}}  + \sum_{k=m+1}^{\infty}\frac{|x_k-\mathbf e(k)|}{2^{k}} <\frac{\varepsilon}{3}+\frac{\varepsilon}{3}+\frac{\varepsilon}{3}=\varepsilon.
$$
We have just proved that the set of end-points of $I_H$ is dense in $I_H$, therefore, $I_H$ is a Lelek fan.
\end{proof}
We conclude this paper by giving {a} proof of our second main result,  Theorem  \ref{mainn}.
\begin{proof}
{  Let $X=I_H^+$.  We have proved in Theorem \ref{Cantor fan plus} that $X$ is a Cantor fan.  Let $f=\sigma_H^+$.  By Theorem \ref{jojmene1}, the mapping $f$  is a transitive continuous surjection.  It follows from Theorem \ref{povezava} that the inverse limit $\varprojlim\{X,f\}$ is homeomorphic to the two-sided Mahavier product $I_H$, which is by Theorem \ref{juuuj} a Lelek fan. Therefore,  the inverse limit $\varprojlim\{X,f\}$ is a Lelek fan. In addition, by Theorem \ref{shifttransitive}, the shift map on $\varprojlim\{X,f\}$ is a transitive homeomorphism.
}
\end{proof}
{  
To our knowledge, there are two known fans that admit a transitive homeomorphism, the Lelek fan and the Cantor fan.  {The Cantor fan may be}  the ``simplest''  non-degenerate one-dimensional continuum  (that is not a simple closed curve) that admits a transitive homeomorphism. Therefore, the following open problems are a good place to finish the paper. 
\begin{problem}
{Are the Cantor fan and the Lelek fan the only fans\footnote{{See \cite[Definition 9, page 5]{banic1} for the definition of a fan.}}  that admit a transitive homeomorphism?}
\end{problem}
\begin{problem}
{Is there ``a simpler'' non-degenerate one-dimensional continuum (that is not a simple closed curve) that admits a transitive homeomorphism}? 
\end{problem}
}

\section{Acknowledgement}
This work is supported in part by the Slovenian Research Agency (research projects J1-4632, BI-HR/23-24-011, BI-US/22-24-086 and BI-US/22-24-094, and research program P1-0285). 
	

\noindent I. Bani\v c\\
              (1) Faculty of Natural Sciences and Mathematics, University of Maribor, Koro\v{s}ka 160, SI-2000 Maribor,
   Slovenia; \\(2) Institute of Mathematics, Physics and Mechanics, Jadranska 19, SI-1000 Ljubljana, 
   Slovenia; \\(3) Andrej Maru\v si\v c Institute, University of Primorska, Muzejski trg 2, SI-6000 Koper,
   Slovenia\\
             {iztok.banic@um.si}           
     
				\-
				
		\noindent G.  Erceg\\
             Faculty of Science, University of Split, Rudera Bo\v skovi\' ca 33, Split,  Croatia\\
{{gorerc@pmfst.hr}       }    

                 	\-
					
  \noindent J.  Kennedy\\
             Department of Mathematics,  Lamar University, 200 Lucas Building, P.O. Box 10047, Beaumont, Texas 77710 USA\\
{{kennedy9905@gmail.com}       }    

	\-
				
		\noindent C.  Mouron\\
             Rhodes College,  2000 North Parkway, Memphis, Tennessee 38112  USA\\\
{{mouronc@rhodes.edu}       }    

                 	\-
				
		\noindent V.  Nall\\
             Department of Mathematics,  University of Richmond, Richmond, Virginia 23173 USA\\
{{vnall@richmond.edu}       }   




\begin{thebibliography}{9}
\bibitem{A} E. ~Akin, {General Topology of Dynamical Systems}, Volume 1, Graduate Studies in Mathematics Series, American Mathematical Society, Providence RI, 1993.
\bibitem{banic1} I.~Bani\v c, G.~Erceg,  J.~Kennedy, The Lelek fan as the inverse limit of intervals with a single set-valued bonding function whose graph is an arc,  { Mediterr. J. Math.  20 (2023) 1--24}
\bibitem{banic2} I.~Bani\v c, G.~Erceg,  J.~Kennedy, A transitive homeomorphism on the Lelek fan,  {to appear in J. Difference Equ. Appl.  (2023) https://doi.org/10.1080/10236198.2023.2208242.}
\bibitem{barge} M.~Barge, J.~Martin, Chaos, periodicity, and snakelike continua. Trans. Amer. Math. Soc. 289 (1985), no. 1, 355--365.
\bibitem{jan} J. ~Boro\' nski, P. ~Minc and S.~ \v Stimac,  On conjugacy between natural extensions of 1-dimensional maps,  Ergod. Th.  Dynam. Sys.   (2022) https://doi.org/10.1017/etds.2022.62.
\bibitem{jan2} J.~ Boro\' nski,  J. ~Kupka, New chaotic planar attractors from smooth zero entropy interval maps, Adv. Difference Equ.  232 (2015) 11 pp.
\bibitem{jan3} J.~ Boro\' nski,  P.~Oprocha, On indecomposability in chaotic attractors. Proc. Amer. Math. Soc. 143 (2015),  3659--3670.
\bibitem{jan4} J.~ Boro\' nski,  P.~Oprocha, On dynamics of the Sierpin\' ski carpet,  C. R. Math. Acad. Sci. Paris 356 (2018) 340--344.
\bibitem{oversteegen} W.~D.~Bula and L.~Overseegen, A Characterization of smooth Cantor Bouquets,  Proc. Amer.Math.Soc. 108 (1990) 529--534.
\bibitem{charatonik} W.~J.~Charatonik, The Lelek fan is unique, Houston J. Math. 15 (1989) 27--34.
\bibitem{cinc} J.~\v Cin\v c, P.~Oprocha, Parametrized family of pseudo-arc attractors: Physical measures and prime end rotations, Proc. London Math. Soc.  125 (2022) 318--357.
\bibitem{handel} M. ~Handel, A pathological area preserving $C^{\infty}$ diffeomorphism of the plane, Proc. Amer.Math.Soc.86 (1982),163--168
\bibitem{he} F. ~He, J. ~Liu, Invariant measures and uniform positive entropy property for inverse limits, Appl. Math. J. Chinese Univ. Ser. B. 14 (1999)  265--272.
\bibitem{HM} L.~C. ~Hoehn and C. ~Mouron, Hierarchies of chaotic maps on continua, Ergodic Theory Dynam. Systems 34 (2014), 1897--1913.
\bibitem{judy} J. ~Kennedy,   A transitive homeomorphism on the pseudoarc which is semiconjugate to the tent map, Trans. Amer. Math. Soc. 326 (1991),  773--793.
\bibitem{ingram}	 W.~T.~Ingram,  An Introduction to Inverse Limits with Set-valued Functions, 	  Springer, New York, 2012.
\bibitem{KS} S.~Kolyada, L.~Snoha,  Topological transitivity, \textit{Scholarpedia}  4 (2):5802 (2009).
\bibitem{lelek} A.~Lelek, On plane dendroids and their end-points in the classical sense, Fund. Math. 49 (1960/1961) 301--319.
\bibitem{li} S. ~Li, Dynamical properties of the shift maps on the inverse limit spaces,  Ergod. Th.  Dynam. Sys.  12 (1992) 95--108.
\bibitem{chris1} V. ~Mart\' nez-de-la-Vega, J.~M. ~Mart\' inez-Montejano, C.~Mouron, Mixing homeomorphisms and indecomposability,  Topology App. \textbf{254} (2019) 50--58.
\bibitem{minc} P. ~Minc and W. ~R.~ R. ~Transue, A Transitive Map on [0,1] Whose Inverse Limit is the Pseudoarc, Proceedings of the American Mathematical Society 111 (1991) 1165--1170. 
\bibitem{chris2} C.~Mouron, Tree-like continua do not admit expansive homeomorphisms. Proceedings of the A.M.S. 130 Nov. 2002, p. 3409-3413.
\bibitem{chris3} C.~Mouron,  Positive entropy homeomorphisms of chainable continua and indecomposable subcontinua, Proc. Amer. Math. Soc. 139 (2011), no. 8, 2783--2791.
\bibitem{chris4} C.~Mouron,  Expansive homeomorphisms and indecomposable subcontinua. Topology Appl. 126 (2002), no. 1-2, 13--28.
\bibitem{chris5} C.~Mouron,  Mixing sets, positive entropy homeomorphisms and non-Suslinean
\bibitem{nadler} S.~B.~Nadler, Continuum theory. An introduction, Marcel Dekker, Inc., New York, 1992.
\bibitem{oprocha} P.~Oprocha,  Lelek fan admits completely scrambled weakly mixing homeomorphism, preprint.
\bibitem{seidler} G.~T.~ Seidler, The topological entropy of homeomorphisms on one-dimensional continua,Proc. Amer. Math. Soc. 4 (1990), 1025--1030.
\end{thebibliography}
\end{document}